\documentclass[a4paper,10pt]{article}
\usepackage{amsmath,amsfonts,amssymb,bbold,bbm,tikz}
\usepackage{algorithm}
\usepackage{algpseudocode}
\usepackage{enumerate}
\usepackage{pgfplots}
\usepackage{pgfplotstable}
\usepackage{booktabs}
\usepackage{geometry}
\usepackage{hyperref}
\usepackage{geometry}
\usepackage{booktabs} 
\usepgfplotslibrary{groupplots}
\pdfoutput=1

\usepackage{cleveref}

\def\bar{\overline}
\def\b{\boldsymbol}

\def\Tr{\mathrm{Tr}}

\newcommand{\norm}[1]{\left\lVert#1\right\rVert}

\newcommand{\upd}[1]{{\color{black}#1}} 

\author{
	Stefano Massei\footnote{Department of Mathematics, University of Pisa. 
		E-mail: stefano.massei@unipi.it. The work of S.M. was partially supported by the INdAM/GNCS project CUP\_E53C22001930001 ``Metodi basati su matrici e tensori strutturati per problemi di algebra lineare di grandi dimensioni".} 
	\and
	Francesco Tudisco\footnote{GSSI  E-mail: francesco.tudisco@gssi.it. The work of F.T. was partially supported by the European Union's Horizon 2020 research and innovation programme under the Marie Sk\l odowska-Curie individual fellowship ``MAGNET'' No 744014.}
}

\title{Optimizing network robustness via Krylov subspaces}

\date{}

\pgfplotsset{compat=1.15}

\begin{document}
\maketitle

\begin{abstract}
We consider the problem of attaining either the maximal increase or reduction of the robustness of a complex network by means of a bounded modification of a subset of the edge weights. We propose two novel strategies combining Krylov subspace approximations with a greedy scheme and \upd{an interior point method employing either the Hessian or its approximation computed via the limited-memory Broyden-Fletcher-Goldfarb-Shanno algorithm (L-BFGS)}. The paper discusses the computational and modeling aspects of our methodology and illustrates the various optimization problems on networks that can be addressed within the proposed framework. 
Finally, in the numerical experiments we compare the performances of our algorithms with state-of-the-art techniques on synthetic and real-world networks.   
\end{abstract}

\section{Introduction}
When studying and analyzing a complex network, one of the main questions is how to identify important nodes and robust connections among them, given the network topology and no other external data. There is a  broad literature on the subject, with many different models and associated algorithms. As a network can be naturally represented by a matrix,  many successful approaches strongly rely on tools from linear algebra and  matrix analysis \cite{chung2006complex, estrada2015first}.

Spectral models, such as eigenvector centrality \cite{vigna2016spectral}, PageRank \cite{gleich2015pagerank}, or resistance distance \cite{luxburg2010getting}, are based on the eigenvalues and eigenvectors of graph matrices and rely on a mutually reinforcing argument, while path-based models, such as Katz centrality \cite{wasserman1994social}, subgraph centrality and total communicability \cite{estrada2008communicability}, use the entries of suitable graph matrix functions and are based on weighted walk counts. 

%
For example,   if $\b x\geq 0$ is the Perron eigenvector of the adjacency matrix $A$ of an undirected graph $G=(V,E)$, then $A\b x = \lambda \b x$ with $\lambda >0$ and hence $x_i$ \upd{is proportional to} $\sum_{j}A_{ij}x_j>0$, for all the nodes $i\in V$. Thus, we can interpret the entries of $\b x$ as an importance score for the nodes of $G$, known as \textit{Bonacich centrality} or \textit{eigenvector centrality} \cite{vigna2016spectral}, where $x_i$, the importance of node  $i$,  is mutually reinforced by the importance of its neighbors. Similarly, if we are given a function of the adjacency matrix %
\begin{equation}\label{eq:pseries}
    f(A) = a_0 I + a_1 A + a_2 A^2 + a_3A^3+\cdots\,\,,
\end{equation}
where the coefficients $a_k$ are  nonnegative, we can interpret the diagonal entries of $f(A)$ as node importances. In fact, in a weighted graph $G$, the weight of a walk from $i$ to $j$ of length $k$ can be defined  as $ A_{u_0u_1}A_{u_1u_2}\cdots A_{u_{k-1}u_k}>0$, where all pairs $u_i u_{i+1}$ are edges and $u_0=i$, $u_k=j$. Thus, the sum of the weights of all the walks of length $k$ from $i$ to $j$ corresponds to $(A^k)_{ij}$ and  the diagonal entry $f(A)_{ii}$ defines the so-called \emph{$f$--centrality} or \emph{subgraph centrality} score of the node $i$ \cite{estrada2010network}, which corresponds to the  weighted sum of all the walks of any length from $i$ and returning to $i$, i.e.\ the subgraphs containing $i$.
Related to the individual centrality of a node are important notions of network robustness and network connectivity, which can be quantified by the summations $\sum_i x_i$,  $\sum_i f(A)_{ii}$ and $\sum_{ij}f(A)_{ij}$, respectively (see e.g.\ \cite{chan2014make,estrada2007statistical}). 
These quantities measure the degree of resiliency of a network in the face of accidental failures or deliberate attacks, modeled as edge modification, removal, or insertion. Both spectral and matrix function-based centrality measures have been successfully used in a variety of settings, including discovering relevant proteins in protein-protein 
interaction networks \cite{estrada2006virtual}, as well as keystone species  in ecological food webs and landscapes \cite{estrada2008using}.

While spectral centralities require the evaluation of one extremal eigenvector and  \upd{can thus} be computed in a relatively cheap way by means of standard sparse numerical eigensolvers, computing the entries of a matrix function can be in general a much more expensive operation, in particular when the matrix is large. This numerical challenge has prompted extensive research work in recent years. Based on Krylov subspace techniques as well as Gauss-Lobatto quadrature formulas, a variety of efficient numerical techniques have been proposed for large-scale sparse networks \cite{alqahtani2018multiple,bellalij2015bounding,fenu2013block,fika2017aitken,hale2008computing,pozza2017stability}.

\upd{Rather than the problem of their efficient evaluation, in this work we focus on the problem of the optimization of matrix function-based node centrality scores.}
Roughly, we look for a ``small'' modification $A+X$ of the current network $A$ that yields the largest centrality increase. 
Here small means that only a limited number of nonzero entries are allowed in $X$ or, in other terms, that we are allowed to modify only a limited number of edges of the graph. Clearly, the resulting optimization task is more complicated than the centrality evaluation problem, as already simple first-order optimization methods would require evaluating both $f(A+X)$ and its Fr\'echet derivative for many different choices of $X$.
Based on recent work on low-rank updates of matrix functions and trace estimators \cite{beck2018update,cortinovis22}, we propose two strategies based on the efficient approximation of $f(A+X)-f(A)$ and the Fr\'echet derivative of $f(A)$ along multiple directions,  to optimize the robustness measure \upd{$\Tr(f(A+X)):=\sum_i f(A+X)_{ii}$}, for both the combinatorial (unweighted) case, in which both $A$ and $A+X$ are binary matrices, and the continuous (weighted) case, in which edge weight tuning is allowed. \upd{Among the most frequently used functions $f$ we mention the exponential function $f(z) = \exp(z)$, which corresponds to the so-called \textit{natural connectivity} \cite{estrada2008communicability}; the hyperbolic sine and consine functions $f(z)=\sinh(z)$, $f(z) = \cosh(z)$, which are often used as a measure of bipartitedness and to define so-called \textit{returnability} \cite{estrada2009returnability}; the resolvent function $f(z) = (1-\alpha z)^{-1}$, which defines the so-called \textit{Katz centrality} \cite{estrada2010network}.} 

The remainder of the paper is structured as follows. In Section~\ref{sec:opt} we introduce the optimization problems that we are going to analyze. Section~\ref{sec:unweighted} describes the greedy  algorithm that we propose in the context of unweighted binary graphs and other techniques that will be used for comparison, see Section~\ref{sec:competitors}. Section~\ref{sec:weight} is dedicated to the gradient method that we propose for weighted graphs. Finally, Section~\ref{sec:experiments} reports numerical experiments concerning optimization problems on both weighted and unweighted graphs.

\subsection{Related work}
Optimizing network robustness or network connectivity  is in general very challenging, due to the combinatorial nature of the problem. A large body of work has focused on spectral-based scores. 
The problem of minimizing the largest eigenvalue (spectral radius) of $A$ by a small number of edge and node removals is considered in  \cite{van2011decreasing,saha2015approximation,zhang2015controlling}. This is shown to be an NP-hard problem which is addressed by a number of heuristics in \cite{van2011decreasing,saha2015approximation} or via a  semidefinite program with polynomial time complexity in \cite{zhang2015controlling}. A similar problem is considered in \cite{tong2012gelling,le2015met}, with the aim of optimizing the network diffusion rate. The works \cite{ghosh2006growing,yu2015friend} studied the problem of maximizing the algebraic connectivity, i.e., the second smallest eigenvalue of the graph Laplacian, and propose both a convex relaxation-based method and a greedy perturbation heuristic, based on the entries of the Fiedler eigenvector of the initial network.   \upd{In \cite{chan2016optimizing} the problem of modifying network edges to reduce external influence is studied.} This is done by controlling the asymptotic consensus value $\b x^T \b a$, where $\b x$ is the eigenvector centrality, i.e.\ the Perron eigenvector of $A$, and $\b a$ is a vector of external user consensus coefficients. The eigenvector centrality $\b x$ is also the subject of \cite{nicosia2012controlling}, where it is observed that, often, modifying a very small subset of edges of a real-world network is enough to drastically change and thus control the eigenvector centrality value of any node in the network. Instead, the Perron eigenvector of the PageRank matrix, so--called PageRank or random walk centrality,  is the subject of \cite{garimella2017reducing}.

Alongside spectral-based coefficients, other network scores have been considered by several authors. For example, \cite{medya2018group} deals with the problem of improving both coverage and betweenness centralities by adding a small set of edges to the network. Greedy algorithms for improving coverage and closeness centralities are proposed in  \cite{d2019coverage} and \cite{crescenzi2016greedily}, respectively.

Centrality optimization problems for indices defined by means of matrix functions are considered for instance in \cite{ghosh2008minimizing,arrigo2016updating,chan2014make}. These works target the optimization of a number of robustness and connectivity coefficients of the network, by modifying, adding, or removing a small subset of edges. In \cite{ghosh2008minimizing}, a semidefinite program-based approach is proposed for the optimization of the total effective resistance, defined as $\Tr(L^+) = \sum_{i}(L^+)_{ii}$, where $L^+$ is the pseudo inverse of \upd{the graph Laplacian} $L$. In \cite{arrigo2016updating,chan2014make}, instead, given a suitable function $f$, a number of heuristics are proposed to efficiently enhance the network natural connectivity, defined as $f^{-1}(\Tr( f(A))/n)$, and the network total communicability \upd{$\b 1^T f(A)\b 1 = \sum_{ij}f(A)_{ij}$}, respectively. Both these two studies show that very good results can be achieved by modifying edges between nodes with high or low centrality values. \upd{The recent work \cite{delacruz22} proposes to measure the sensitivity of the network communicability, to the addition or removal of certain edges, by looking at the derivatives of $\b 1^T f(A)\b 1$. The latter quantities called \emph{total network sensitivities}, are defined in terms of evaluations of the Fr\'echet derivative of $f(A)$.  The preprint by Schweitzer \cite{schweitzer23}, which appeared in parallel to the first version of this document, introduces an efficient technique that is able to compute all the total network sensitivities by means of a single evaluation of the Frech\'et derivative of $f(A)$ in the rank one direction $\b 1 \b 1^T$. An analogous technology is applicable for computing the derivatives of the network's natural connectivity.}

Building on top of this body of work,  we focus here on the  optimal modification of the network's natural $f$-connectivity. In the sequel, we formalize the problem and the algorithms we propose.

\section{Optimizing the natural connectivity}\label{sec:opt}
Networks strongly rely on their robustness, i.e., the ability to maintain a high degree of connectivity when a portion of the network's structure is damaged or simply altered. An intuitive notion of graph robustness can be expressed in terms of the redundancy of routes between vertices. If we consider a source vertex and a termination vertex, there may be several paths between them. When one path fails, the two vertices can still communicate via other alternative routes. \upd{Hence, the robustness of the network grows with the number of available alternative routes.}
Thus,  an ideal measure of robustness for a network would be the degree of  redundancy of alternative paths, i.e. the number of alternative routes of different lengths for all pairs of vertices. However, this number is very difficult to compute. 

An alternative definition of robustness, which is usually called ``natural connectivity'', counts instead the number of closed walks of any length. 
Let $G=(V,E)$ be an undirected, possibly weighted graph with $V=\{1,\dots,n\}$ and entry wise nonnegative symmetric adjacency matrix $A\geq 0$, such that $A_{ij}>0$ if and only if $ij\in E$. As the number of closed walks of length $k$ from $i$ to itself coincides with the $i$-th diagonal entry of the $k$-th power of the adjacency matrix, we can quantify the natural connectivity  by looking at 
\[
\ln \left(\frac 1 n \sum_{i=1}^n \sum_{k=1}^\infty \frac {\left(A^k\right)_{ii}}{k!}\right) = \ln \Big(\frac 1 n \Tr(\exp(A))\Big) = \ln\Big( \frac 1 n \sum_{i=1}^n e^{\lambda_i}\Big)
\]
where $\lambda_1\leq \dots\leq\lambda_n$ are the eigenvalues of the adjacency matrix $A$. 
The scaling factor $1/k!$ is required here in order to have a convergent series \upd{and to discount the importance of long walks with respect to short ones.} The logarithm and the scaling factor $1/n$ are used to avoid very large numbers as they yield an ``average'' of the eigenvalues of the adjacency matrix. More in general, we can consider the natural $f$-connectivity ($f$-connectivity, in short) 
as the generalized $f$-mean of the eigenvalues of $A$
\[
\vartheta(A)  = f^{-1}\Big(\frac 1 n \Tr(f(A))\Big)= f^{-1}\Big( \frac 1 n \sum_{i=1}^n f(\lambda_i)\Big),
\]
where $f$ is a real-valued, increasing, and analytic function on a set containing the spectrum of~$A$. 

As $f$ is increasing, it is not difficult to realize that $\vartheta(A)$ itself changes monotonically with the edges of the graph, that is, $\vartheta(A)$ grows if edges are added, and decreases if they are removed. 
In the following, we  assume we are given a \textit{budget} $k$ representing the number of edges, or the cumulative edges' weight,  that can be either removed or added to the graph. Thus, we consider  the \upd{optimization} problem of using the given budget to either reduce or increase $\vartheta(A)$ the most. 

In matrix terms we can formulate the corresponding optimization problem as follows. Assume we are given the initial graph with adjacency matrix $A$.
We want to find a modification $X$ of the network edges $A$ that either maximizes or minimizes the function $\vartheta(A+X)$, subject to suitable  constraints on $X$ which account for the budget and for whether we are removing, adding or modifying the weight of the edges, as detailed next. The constraints on $X$ also depend on whether we are considering weighted or unweighted (binary) networks. To summarize we consider the following three classes of optimization problems.

\subsubsection*{Edge downgrading}
Let us assume that we are given a \upd{positive budget $k$} and we want to remove or diminish the weight of the edges that yield the greatest decrease in $f$-connectivity. Given the graph $G=(V,E)$, we then consider the set of admissible modifications 
$$
\Omega_k(E) = \Big\{X : \textstyle{\sum_{ij}|X_{ij}|\leq k}, \text{ $X=X^T$,  $X_{ij} = 0$, for $ij \notin E$}\Big\}.
$$
The downgrading problem for  unweighted graphs, more often referred to as \textit{edge breaking} problem \cite{chan2014make},   is:
\begin{equation}\label{eq:P-removal}\tag{DG}
    \min          \,  \vartheta(A+X) \;\;
     \text{s.t.} \;\;   X\in \Omega_k(E) \text{ and } X_{ij}\in \{-1,0\}  
\end{equation}
while for weighted graphs the second constraint is  replaced by $-A_{ij}\leq X_{ij}\leq 0$, i.e.
\begin{equation}\label{eq:P-downgrading}\tag{DG'}
  \min          \,  \vartheta(A+X) \;\;
     \text{s.t.} \;\;   X\in \Omega_k(E) \text{ and }   -A_{ij}\leq X_{ij}\leq 0. 
\end{equation}

\subsubsection*{Edge addition}
In this setting, we consider the situation where new edges may be introduced in order to increase the $f$-connectivity of the network. In this case, given a budget $k$, the set of admissible modifications takes the form
\[
 \Omega_k(\bar E) = \Big\{X : \textstyle{\sum_{ij}X_{ij}\leq k}, \text{ $X=X^T$,  $X_{ij} = 0$, for $ij \in E$}\Big\}.
\]
For unweighted graphs, we obtain the following optimal edge addition problem 
\begin{equation}\label{eq:P-addition-binary}\tag{AD}
  \max          \,  \vartheta(A+X) \;\;
     \text{s.t.} \;\;   X\in  \Omega_k(\bar E)   \text{ and } X_{ij}\in \{0,1\}.
\end{equation}
To avoid trivial solutions, where all the budget is spent on a single most important edge, when dealing with  weighted networks, we further assume we are given a set of maximum weight values $U_{ij}$ that we are allowed to spend on each edge:
\begin{equation}\label{eq:P-addition-weighted}\tag{AD'}
  \max          \,  \vartheta(A+X) \;\;
     \text{s.t.} \;\;   X\in \Omega_k(\bar E) \;\; 0\leq X_{ij}\leq U_{ij}.
\end{equation}

\subsubsection*{Edge tuning}
Finally, in the third problem, we are given the budget $k$ and a weighted graph $G$, and we look for a modification of the edge weights of a limited set $F\subseteq E$ of the existing edges in order to obtain the  largest increase in $f$-connectivity. We will also consider the case where $F$ includes both existing and non existing edges, to address the scenario where the creation of new links is also allowed; we call this slightly modified problem \emph{edge rewiring}. As for \eqref{eq:P-addition-weighted}, we  assume a set of maximum weight values $U_{ij}$ is given,  to avoid trivial solutions: 
\begin{equation}\label{eq:P-tuning}\tag{TU}
  \max          \,  \vartheta(A+X) \;\;
     \text{s.t.} \;\;  \text{$X\in \Omega_k(F)$ and }   -A_{ij}\leq X_{ij}\leq U_{ij} .
\end{equation}

\subsection{Algorithmic set-up}
Before moving on to the proposed algorithmic techniques, we make several preliminary remarks.

First, we note that since $f$ in the definition of $\vartheta$ is an increasing function, \upd{then so is $f^{-1}$. Additionally, note that } minimizing (resp.\ maximizing) the natural $f$-connectivity is equivalent to minimizing (resp.\ maximizing) the trace variation 
\[
\varphi_A(X) := \Tr(f(A+X))-\Tr(f(A))\, ,
\]
with respect to $X$.

\upd{Secondly,} we observe that the dimensions of the constraint sets that involve all the existing (or non-existing) edges in the graph are usually very large, already for  graphs of moderate size. For this reason, in the rest of the paper we further restrict the optimization problems above to a subset of the edges (or non-existing edges)  $F$ whose  elements are cleverly selected and whose size is kept under control. 

The selection of a suitable $F$ may depend on the problem at hand, and we will call this procedure ``the search space selection'', which will be discussed \upd{case-by-case} in Sections \ref{sec:search} and \ref{sec:exp-weight}. 
Note that, in real applications, further constraints on the set of modifiable edges (or non-existing edges) may be imposed by the application set-up: for example, one may have only access to a certain part of the network (\upd{as in the case} of a street network where most of the  roads may not be modifiable). \upd{This} additional problem-based constraint can be \upd{imposed} by straightforward modifications of the above optimization problems.

\section{Edge downgrading and addition for unweighted graphs}\label{sec:unweighted}
In this section we propose some heuristic greedy procedures  for addressing the optimization problems \eqref{eq:P-removal} and  \eqref{eq:P-addition-binary}. We begin by describing the general greedy template that is behind our method and other algorithms proposed in the literature. Throughout the discussion we assume to have a budget of $k$ edges.
\subsection{The greedy paradigm}
 The most intuitive greedy strategy for problem \eqref{eq:P-removal} (resp.\ \eqref{eq:P-addition-binary}) consists of sequentially removing (resp.\ adding) the edge that attains the largest reduction (resp.\ increase) of $\varphi_A$ until $k$ deletions (resp.\ additions) are performed. Usually, the identification of the $j$th edge to be either added or removed is made by evaluating or approximating the variation of $\varphi_A$ on  a large number of candidate edges. 
Even in the case of an exhaustive search of candidates over the whole edge set (or the whole set of missing edges, in the case of \eqref{eq:P-addition-binary}), this greedy procedure is guaranteed to return the optimal solution only for $k=1$; on the other hand, when $k>1$, we expect that the selected set of $k$ edges provides a significant modification of $\varphi_A$.

When dealing with medium to large networks, the implementation of this greedy procedure  poses two major computational issues: 
\begin{itemize}
    \item[$(i)$] the large number of edges in the search space to be processed in each step, and
    \item[$(ii)$] the cost of evaluating (or approximating) the cost function $\varphi_A$.
\end{itemize}
 Concerning $(i)$, we remark that, when the graph is sparse, an exhaustive search would require considering 
 $\mathcal O(n)$ edges for problem \eqref{eq:P-removal} and $\mathcal O(n^2)$  edges for problem \eqref{eq:P-addition-binary}. When such sets have large sizes, this step can be prohibitively expensive. This is circumvented by restricting the search space for the $j$th edge  to an appropriate subset $F_j$ of moderate size.
In the case of \eqref{eq:P-removal}, $F_j\subseteq E$, while for \eqref{eq:P-addition-binary} $F_j\subseteq V\times V\setminus E$.

 Similarly, task $(ii)$ involves $f(A)$ and $f(A+X)$ but \upd{cannot} be addressed by directly forming these matrix functions as \upd{$f(A)$ is dense almost always, even if $A$ is sparse, and computing $f(A)$ directly would require $\mathcal O(n^3)$ operations.} 
 Even for small to medium-size matrices, as $X$ changes at each greedy step, computing $f(A+X)$ each time would be prohibitively expensive. Efficient greedy methods make use of techniques that approximate the variation  $\varphi_A$ with a reduced computational cost.

 In Algorithm~\ref{alg:greedy} we present a general scheme for the above greedy strategy, for the case of \eqref{eq:P-removal}. The analogous algorithm for \eqref{eq:P-addition-binary} is obtained with straightforward modifications at lines \ref{line:rank2} and \ref{line:crit}, by changing the sign of the rank-2 update and  reversing the inequality for $\delta_{\sf opt}$, which has to be  initially set to $-\infty$ \upd{at line \ref{step:delta-opt}}.
Then, in the next two subsections, we will present our proposed strategy for addressing the two points $(i)$ and $(ii)$ above. In particular, we propose the use of a Krylov subspace-based approach for the approximation of the variation $\varphi_A$, which will guarantee an accurate approximation with a computational cost of $\mathcal O(n)$, \upd{as detailed in subsection \ref{sec:update}}.

 \begin{algorithm}[t] \small
	\caption{Template of a greedy method for \eqref{eq:P-removal}}\label{alg:greedy}
	\begin{algorithmic}[1]
	\Procedure{greedy\_downgrade}{$A$, $k$}
	\State Set $\Delta A=0$
		\For{$j = 1, \ldots, k$}
		\State $X_{\mathsf{opt}} \gets 0$, $\delta_{\mathsf{opt}}\gets +\infty$ \label{step:delta-opt}
		\State Select $F_j$
		\For{$(s,t)\in F_j$}\label{line:Ej}
		\State \label{line:rank2} $X\gets -({\bf 1}_s {\bf 1}_t^T+{\bf 1}_t{\bf 1}_s^T)$\Comment{rank $2$ modification that deletes $(s,t)$}
		\State Compute $\varphi_A(X) = \Tr(f(A+X))-\Tr(f(A))$\label{line:crit}
		\If{$\varphi_A(X)\leq \delta_{\mathsf{opt}}$}
		\State $\delta_{\mathsf{opt}}\gets \varphi_A(X)$
		\State $X_{\mathsf{opt}}\gets X$
		\EndIf
		\EndFor
				\State $\Delta A \gets \Delta A + X_{\mathsf{opt}} $
	\State \label{line:adj_update}$A\gets A+ X_{\mathsf{opt}}$
		\EndFor
		\State \Return $\delta_{\mathsf{opt}},\Delta A$
		\EndProcedure
	\end{algorithmic}
\end{algorithm}

\subsection{Selection of the search spaces}\label{sec:search}
The strategy for selecting the sets $F_j$ has to ensure a feasible size of the search space and that the most meaningful edges are considered. Intuitively, the second requirement is the trickiest as, due to the combinatorial nature of \eqref{eq:P-removal} and \eqref{eq:P-addition-binary}, only an exhaustive search space can guarantee it. The latter choice might be computationally viable for problem \eqref{eq:P-removal} where each $F_j$ has at most $\mathcal O(n)$ edges, assuming the initial graph is sparse. If $n$ is moderate and the cost of evaluating $\varphi_A(X)$ 
is at most linear on $n$, then we consider the following search spaces
\begin{equation}\label{strat:exhaustive}\tag{\text{$\text{S}_{\text{DG}}^{\text{full}}$}}
\begin{cases}
F_1=E\\
F_{j+1}=E\setminus \mathsf{Chosen}(j)
\end{cases},
\end{equation}
with $\mathsf{Chosen}(j):=\{\text{edges selected in the first $j$-th steps of Algorithm~\ref{alg:greedy}}\}$.

When strategy \eqref{strat:exhaustive} is too expensive, an alternative is to define a ranking on the set of edges to heuristically identify the most important ones. Here we propose to rank the edges on the basis of the eigenvector centrality scores of the nodes they connect, as these scores for the nodes are cheap to evaluate for sparse graphs.   More specifically, given two edges $(v_1,v_2)$ and $(v_3,v_4)$, we consider the following two rankings $\leq_1$ and $\leq_2$ on $V\times V$:
\begin{align*}
    (v_1,v_2)\leq_1(v_3,v_4)\quad&\Longleftrightarrow\quad \mathsf{eigc}(v_1)\cdot\mathsf{eigc}(v_2)\leq \mathsf{eigc}(v_3)\cdot \mathsf{eigc}(v_4),\\
    (v_1,v_2)\leq_2(v_3,v_4)\quad&\Longleftrightarrow\quad \begin{cases}\min\{\mathsf{eigc}(v_1),\mathsf{eigc}(v_2)\}< \min\{\mathsf{eigc}(v_3),\mathsf{eigc}(v_4)\}\\
    \qquad\qquad\qquad\qquad\qquad\text{ or}\\
    \min\{\mathsf{eigc}(v_1),\mathsf{eigc}(v_2)\} = \min\{\mathsf{eigc}(v_3),\mathsf{eigc}(v_4)\}\\
    \max\{\mathsf{eigc}(v_1),\mathsf{eigc}(v_2)\}\leq \max\{\mathsf{eigc}(v_3),\mathsf{eigc}(v_4)\}
    \end{cases}.
\end{align*}
where  $\mathsf{eigc}(v)$ denotes the eigenvector centrality of node $v\in V$, i.e.\ the $v$-th entry $x_v$ of the  Perron eigenvector $\b x$ of the adjacency matrix. The ordering $\leq_1$ is a standard way of inferring centralities for edges from the node scores \cite{arrigo2016updating, tudisco2021node}. However, we note that $\leq_1$ may still assign large importance to edges that connect \upd{a node with small centrality with another having a large centrality}; this is prevented by $\leq_2$ which thresholds the edge score by the smallest node centrality involved.  We observe that $\leq_2$ works better in practice, as shown in  the numerical experiments   in Section \ref{sec:experiments}.

Finally, given a subset of edges $F\subseteq V\times V$ and a positive integer $q$, we denote with $\left[F\right]_q^{\leq_i}$ the subset of $F$ made by its largest $q$ elements according to $\leq_i$, $i=1,2$. The following selection strategies maintain a search space of size $q$ at each step of Algorithm~\ref{alg:greedy}: \\
\begin{minipage}{.44\textwidth}
\begin{align}
    \begin{cases}
    F_1= \left[E\right]_q^{\leq_1}\\
    F_{j+1} =\left[E\right]_{q+j}^{\leq_1}\setminus \mathsf{Chosen}(j)
    \end{cases}\label{strat:mult-del}\tag{\text{$\text{S}_{\text{DG}}^{1}$}}\\
        \begin{cases}
    F_1= \left[E\right]_q^{\leq_2}\\
    F_{j+1} =\left[E\right]_{q+j}^{\leq_2}\setminus \mathsf{Chosen}(j)
    \end{cases}\label{strat:min-del}\tag{\text{$\text{S}_{\text{DG}}^{2}$}}
\end{align}
\end{minipage}
\begin{minipage}{.5\textwidth}
\begin{align}
    \begin{cases}
    F_1= \left[V\times V\setminus E\right]_q^{\leq_1}\\
    F_{j+1} =\left[V\times V\setminus E\right]_{q+j}^{\leq_1}\setminus \mathsf{Chosen}(j)
    \end{cases}\label{strat:mult-add}\tag{\text{$\text{S}_{\text{AD}}^{1}$}}\\
        \begin{cases}
    F_1= \left[V\times V\setminus E\right]_q^{\leq_2}\\
    F_{j+1} =\left[V\times V\setminus E\right]_{q+j}^{\leq_2}\setminus \mathsf{Chosen}(j)\label{strat:min-add}\tag{\text{$\text{S}_{\text{AD}}^{2}$}}
    \end{cases}
\end{align}
\end{minipage}\\[.3em]
where we have used the subscripts DG and AD to emphasize that the corresponding strategy is meant for  problem \eqref{eq:P-removal} and \eqref{eq:P-addition-binary}, respectively.

Finally, we describe an additional selection strategy for \eqref{eq:P-addition-binary} proposed in \cite{chan2014make}, a method we will use as benchmark for comparison in our experiments. Let $d$ be the  maximum node degree  of the graph and denote by $V_d\subseteq V$ the set  of $d$ nodes of largest degrees. Then, the selection strategy uses the missing edges contained in $V_d\times V_d$. This is formally expressed with the following equation:
\begin{equation}\label{strat:degree}\tag{\text{$\text{S}_{\text{AD}}^{3}$}}
\begin{cases}
F_1=V_d\times V_d\setminus E\\
F_{j+1}=V_d\times V_d\setminus \{E\cup\mathsf{Chosen}(j)\}
\end{cases}.
\end{equation}
Note that, strategy \eqref{strat:degree} only ensures that the search space has cardinality bounded from above by $d^2$; this might be a very weak property for certain graph topologies, as $|F_{j}|$ can be very small.

\subsection{Updating the trace of $f(A)$}\label{sec:update}
The main computational efforts of \Cref{alg:greedy} 
come from evaluating $\varphi_A(X)$ 
at line~\ref{line:crit}. Note that the matrix $X$ at that step of the algorithm is symmetric and has rank 2. Leveraging this key rank property, we can devise a method of cost $\mathcal O(n)$ for computing the variation $\varphi_A(X)$, 
based on the Krylov subspace method in \cite{beck2018update}. 
We start by describing in Section \ref{sec:trace_update}  the proposed Krylov method; then, in Section \ref{sec:miobi} we report another approximation of $\varphi_A$ that has been previously used in the literature and that will be used as a baseline for comparison later. 
\subsubsection{A Krylov projection method}\label{sec:trace_update}
 Let $A$ be a symmetric adjacency matrix, $X$ a symmetric low-rank modification and $f(z)$ a scalar function.
In \cite{beck2018update} it has been proved that, under mild assumptions, the matrix $\Delta f:=f(A+X)-f(A)$ is  of low numerical rank and its approximation can be performed by means of Krylov subspaces. We will see that, with some minor modifications, this also allows to cheaply approximate  $\Tr(\Delta f)=\Tr(f(A+X))-\Tr(f(A))$.
 
 Let us assume $X=\b U_X B_X\b U_X^*$ with $\b U_X\in\mathbb R^{n\times s}$, $B_X=B_X^*\in\mathbb R^{s\times s}$ and denote by 
 $\mathcal{K}_m(A, \b U_X)$ the $m$-th order \emph{Krylov subspace} generated by $A$ and the (block) vector $\b U_X$: 
 $$ \mathcal{K}_m(A, \b U_X):=\text{Span}\{\b U_X, A\b U_X, \dots, A^{m-1}\b U_X\},$$
 \upd{where $\text{Span}$ indicates the column span.}
 If $m$ steps of the Arnoldi process on $A$ and $\b U_X$ can be carried out without breakdowns, then it returns an orthonormal basis $\b{\mathcal U}_m=[\b U_1 |\dots|\b U_m]\in\mathbb R^{n\times ms}$  of $\mathcal{K}_m(A, \b U_X)$ which verifies the following block Arnoldi relation \cite{Gutknecht2007}:
\begin{equation} \label{eq:arnoldirelations}
 A \b{\mathcal U}_m = \b{\mathcal U}_m \mathcal H_m + \b U_{m+1} H_{m+1,m} \b E_m^T,
\end{equation}
with a $ms\times ms$ block tridiagonal matrix $\mathcal H_m$, a $s\times s$ matrix $H_{m+1,m}$,
and $\b E_m^T = [0 | \cdots | 0 | I_s]\in\mathbb R^{s\times ms}$, where $I_s$ denotes the $s\times s$ identity matrix. An approximation of $\Delta f$ is given by 
\begin{equation}\label{eq:deltaf-approx}
    \Delta f\approx \Delta_m f := \b{\mathcal  U}_m [f(\mathcal H_m+\b W_m B_X\b W_m^*)-f(\mathcal H_m)]\b{\mathcal  U}_m^*,
\end{equation}
where $\b W_m :=\b{\mathcal  U}_m^*\b U_X\in\mathbb R^{ms\times s}$. The algorithm proposed in \cite{beck2018update}, reported in Algorithm~\ref{alg:fun-update}, builds --- incrementally in $m$ --- the Arnoldi relations \eqref{eq:arnoldirelations} and their corresponding quantities $\Delta_m f $. We remark that the matrix $\Delta_m f $ is kept in the factored form $\Delta_m f = \b{\mathcal  U}_m \widetilde{\Delta}_m f\ \b{\mathcal  U}_m^*$ where $\widetilde{\Delta}_m f :=f(\mathcal H_m+\b W_m B_X \b W_m^*)-f(\mathcal H_m)\in\mathbb R^{ms\times ms}$.  The method stops when the heuristic stopping criterion $$\norm{\Delta_m f -\Delta_{m-\ell}f}_2=\norm{\widetilde{\Delta}_m f -\begin{bmatrix}\widetilde{\Delta}_{m-\ell}f& 0 \\ 0 &0\end{bmatrix}}_2\leq \epsilon$$ is satisfied for a prescribed tolerance $\epsilon$ and a positive integer $\ell$; in our implementation we set $\ell=2$. \upd{We emphasize that this is just one (arguably, the simplest) of a variety of possible choices for the  stopping criterion. Alternative and more accurate methods for computing error estimates of block Arnoldi methods for matrix functions can be used, as discussed for example in \cite{chen2022error,frommer2017block}.} 

Concerning the approximation error, the method is exact when $f(z)$ is a low degree polynomial; more precisely, $\Delta f=\Delta_m f$ when $f\in\mathcal P_{m-1}$, where $\mathcal P_{m-1}$ denotes the set of polynomials of degree at most $m-1$. For a more general $f$, the error norm is linked to the best polynomial approximation of $f$ on a set $\Pi$ containing the convex hull of the spectrum of $A$ and $A+X$ \cite[Theorem 4.1]{beck2018update}.

We remark that, if the goal is to approximate $\Tr(\Delta f)$, then we can avoid the evaluation of matrix functions at all. Indeed, for computing $\Tr(\Delta_m f)= \Tr(f(\mathcal H_m+\b W_m B_X \b W_m^*))-\Tr(f(\mathcal H_m))$ it is sufficient to retrieve the eigenvalues of the small symmetric matrices $\mathcal H_m$ and $\mathcal H_m+\b W_m B_X \b W_m^*$, and then apply the function $f$ to them. \upd{Since only the approximate eigenvalues are needed here, we replace the Arnoldi method with the Lanczos method for computing the projected matrices.} Moreover, a tighter approximation bound is obtained for this particular case, namely \cite[Theorem 3]{cortinovis22}:
\[
|\Tr(\Delta f) - \Tr(\Delta_m f)|\leq 4n\min_{p\in\mathcal P_{2m}}\max_{z\in \Pi} |f(z)-p(z)|.
\]
We report the pseudocode of the procedure for approximating the variation $\Tr(f(A+X))-\Tr(f(A))$ in Algorithm~\ref{alg:trace-fun-update}. 

\upd{Under the assumptions that matrix-vector products with the matrix $A$ cost $\mathcal O(n)$, that the rank of $X$ is $r$, and that $\mathrm{it}$ iterations of the Arnoldi method have been executed before detecting convergence, the cost of Algorithm~\ref{alg:fun-update} is $\mathcal O(n r^2\mathrm{it}^2+r^3\mathrm{it}^4)$. The term of complexity $\mathcal O(n r^2\mathrm{it}^2)$ comes from the full re-orthogonalization applied in the Arnoldi procedure; moreover, computing $\widetilde{\Delta}_m$ requires the evaluation of two functions of
	$rm \times rm$ symmetric matrices, which typically needs $\mathcal O(r^3m^3)$, and this yields the term of complexity $\mathcal O(r^3\mathrm{it}^4)$. An analogous analysis applies to Algorithm~\ref{alg:trace-fun-update} that is of complexity $\mathcal O(nr\ \mathrm{it}+r^3\mathrm{it}^4)$; the major difference with Algorithm~\ref{alg:fun-update}, is that at, each iteration, the Lanczos method only orthogonalizes with respect to the last two block vectors of the orthonormal basis, and that the eigenvalues of two $rm\times rm$ symmetric matrices are computed in place of their matrix functions. Note that, when calling Algorithm~\ref{alg:trace-fun-update} from Algorithm~\ref{alg:greedy} we always have $r=2$.}
 \begin{algorithm} \small
	\caption{Low-rank approximation of $f(A+X)-f(A)$\label{alg:fun-update}}
	\begin{algorithmic}[1]
		\Procedure{fun\_update}{$A$, $\b U_X$, $B_X$, $f$, $\ell$, $\epsilon$}
		\For{$m = 1, \ldots, m_{\max}$}
		\State \begin{minipage}{0.8\textwidth} Compute (incrementally) the Arnoldi relation for $\mathcal K_{m}(A,\b U_X)$ \upd{by means of the Arnoldi method};  store 
			$\b{\mathcal U}_{m}=[\b U_1 | \dots | \b U_{m}]$ and $\mathcal H_m$
		\end{minipage}
		\State $\b W_m \gets \b{\mathcal U}_{m}^*\b U_X$
		\State $\widetilde{\Delta}_mf\gets f(\mathcal H_m+\b W_m B_X \b W_m^*)- f(\mathcal H_m)$
		\If{$m> \ell$ and $\norm{\widetilde{\Delta}_m f -\begin{bmatrix}\widetilde{\Delta}_{m-\ell}f& 0 \\ 0 &0\end{bmatrix}}_2\leq \epsilon$}
		\State \textbf{break}
		\EndIf
		\EndFor
		\State \Return $\b {\mathcal U}_m, \widetilde{\Delta}_mf$
		\EndProcedure
	\end{algorithmic}
\end{algorithm}
 \begin{algorithm} \small
	\caption{Approximation of $\Tr(f(A+X)-f(A))$\label{alg:trace-fun-update}}
	\begin{algorithmic}[1]
	\Procedure{trace\_fun\_update}{$A$, $\b U_X$, $B_X$, $f$, $\ell$, $\epsilon$}
		\For{$m = 1, \ldots, m_{\max}$}
		\State \begin{minipage}{0.8\textwidth} Compute (incrementally) the Arnoldi relation for $\mathcal K_{m}(A,\b U_X)$ \upd{by means of the Lanczos method};  store 
		 $\b{\mathcal U}_{m}=[\b U_1 | \dots | \b U_{m}]$ and $\mathcal H_m$
		 \end{minipage}
		 \State $\b W_m \gets \b{\mathcal U}_{m}^*\b U_X$
		 \State Compute the eigenvalues $\widetilde \lambda_j$ of $\mathcal H_m+\b W_m B_X \b W_m^*$
		 \State Compute the eigenvalues $\lambda_j$ of $\mathcal H_m$
		 \State $\Delta_m\lambda \gets \sum_j f(\widetilde \lambda_j)-f(\lambda_j)$
		 \If{$m> \ell$ and $|\Delta_m\lambda -\Delta_{m-\ell}\lambda|< \epsilon$}
		 \State \textbf{break}
		 \EndIf
		\EndFor
		\State \Return $\Delta_m\lambda$
		\EndProcedure
	\end{algorithmic}
\end{algorithm}
\subsubsection{Approximation via eigendecomposition update}\label{sec:miobi}
The algorithm \emph{make it or break it} (MIOBI) proposed in \cite{chan2014make} approximates the difference of traces by means of a first-order approximation of the largest eigenpairs of $A+X$. More specifically, given a positive integer $h$, the procedure starts by computing the eigenpairs $(\lambda_1,\b u_1),\dots,(\lambda_h, \b u_h)$ of $A$, corresponding to the $h$ eigenvalues of largest magnitudes. For each $X$, the authors of \cite{chan2014make} observe that the dominant $h$ eigenpairs $\widehat \lambda_j$, $\widehat{\b  u_j}$ of $A+X$ can be written as 
\begin{equation*}\label{eq:update-formula}
\begin{array}{l}
   \widehat\lambda_j= \widetilde
    \lambda_j + \mathcal O(\|X\|^2)\\
   \widehat{\b u}_j=\widetilde{\b u}_j+ \mathcal O(\|X\|^2) 
\end{array}
\quad \text{with} \quad \begin{array}{l}
   \widetilde \lambda_j= \lambda_j + \b u_j^*X\b u_j,\\
   \widetilde {\b u}_j=\b u_j +\sum_{i=1,i\neq j}^h\frac{\b u_i^*X\b u_j}{\lambda_i-\lambda_j}\b u_i.
\end{array}
\end{equation*}
Thus, it is proposed  to consider  the pairs $(\widetilde{\lambda}_j, \widetilde{\b u}_j)$ as approximations of $(\widehat{\lambda}_j, \widehat{\b u}_j)$, i.e., to neglect the high-order terms $\mathcal O(\|X^2\|)$. \upd{This approach is particularly useful when $X$ is a perturbation with small norm.} The resulting procedure is of the same form as Algorithm~\ref{alg:greedy}, with two main modifications: at line \ref{line:crit} the  formula $\sum_{j=1}^h f(\widetilde \lambda_j)-f(\lambda_j)$ is used to approximate the trace update $\Tr(f(A+X))-\Tr(f(A))$; then at line \ref{line:adj_update} both  formulas for $\widetilde{\lambda}_j, \widetilde{\b u}_j$ are used to approximate the $h$ dominant eigenpairs of $A+X_{\sf opt}$. Overall, this yields an algorithm with an iteration cost of $\mathcal O(|F_j|h + nh^2)$.
\subsection{Algorithms for edge downgrading and edge addition}\label{sec:competitors}
We are now ready to formally introduce the methods that we propose for solving \eqref{eq:P-removal},\eqref{eq:P-addition-binary}:
\begin{description}
\item[\normalfont\textsc{greedy\_krylov\_break:}] Algorithm~\ref{alg:greedy} combined with \textsc{trace\_fun\_update} for evaluating the difference of traces at line~\ref{line:crit} and using the strategy \eqref{strat:min-del} for selecting the sets $F_j$.
\item[\normalfont\textsc{greedy\_krylov\_make:}] Algorithm~\ref{alg:greedy} combined with \textsc{trace\_fun\_update} for evaluating the difference of traces at line~\ref{line:crit} and using the strategy \eqref{strat:min-add} for selecting the sets $F_j$.
\end{description}

Moreover, to provide a comparison with the performance of state-of-the-art greedy schemes, we consider the following methods:

\begin{description}\item[\normalfont\textsc{miobi:}] Greedy method proposed in \cite{chan2014make} that uses \eqref{eq:update-formula} for evaluating the difference of traces and the selection strategies \eqref{strat:exhaustive} and \eqref{strat:degree} for \eqref{eq:P-removal} and \eqref{eq:P-addition-binary}, respectively.
\item[\normalfont\textsc{eigenv:}] Method proposed in \cite{arrigo2016updating} that consists in deleting or adding the $k$ edges with the largest eigenvector centrality scores --- with respect to $\le_1$ --- in $E$ and $V\times V\setminus E$, respectively. \upd{The dominant part of its cost is given by the computation of the dominant eigenvector of the adjacency matrix; in our implementation this is done by means of the Matlab function \texttt{eigs}.}
\end{description}

\section{Edge downgrading, addition, and tuning for weighted graphs}\label{sec:weight}
When considering the solution of \eqref{eq:P-downgrading}, \eqref{eq:P-addition-weighted}, and \eqref{eq:P-tuning}, one needs to deal with a constrained continuous optimization problem involving the objective function $\varphi_A(X)$. Similarly to what has been done for unweighted graphs in Section~\ref{sec:unweighted}, we keep  the size of the problem \upd{under control} by imposing that we are allowed to modify only a subset $F$ of the edges (or the missing edges), with  cardinality $n_F$. With this constraint, we have that $\varphi_A$ can be seen as a function of $n_F$ variables $\varphi_A:\mathbb R^{n_F}\rightarrow \mathbb R$, which correspond to the variation of the weights of the edges in $F$. In particular, the matrix $X$ has rank bounded by $2n_F$; i.e., to efficiently evaluate $\varphi_A(X)$ we can rely on Algorithm~\ref{alg:trace-fun-update}, as far as $2n_F\ll n$.

We perform the efficient optimization of $\varphi_A$  via \upd{two tailored implementations of an Interior-Point method. The first one, approximates the Hessian of the objective function by means of the Limited-memory BFGS algorithm (L-BFGS),  which iteratively updates the approximation  via rank-2 corrections and only requires the evaluation of the objective function and its gradient. The second one, approximates the true Hessian by means of a Krylov approach. Note that the second approach involves the computation of the second derivatives while the first approach does not. The evaluations of $\varphi_A(X)$ are  computed by means of
Algorithm~\ref{alg:trace-fun-update} as in the discrete setting. The gradient and Hessian computations, instead, require additional analysis as they can be prohibitively expensive if done in a naive way.  We devote  the remainder of this section to briefly review the L-BFGS algorithm and to the description of numerical methods, presented in Algorithm~\ref{alg:gradient} and Algorithm~\ref{alg:hessian}, to efficiently evaluate the gradient and the Hessian of $\varphi_A$. The ultimate procedures obtained by combining the Interior-Point method with Algorithm~\ref{alg:trace-fun-update} for the objective function evaluation, Algorithm~\ref{alg:gradient} for the gradient, and either L-BFGS or Algorithm~\ref{alg:hessian} for the Hessian, are denoted by \textsc{krylov\_lbfgs} and \textsc{krylov\_hessian}, respectively. Our implementation of the Interior-Point method relies on the Matlab function \texttt{fmincon} that allows us  to specify handle functions for the evaluation of the objective function, the gradient, and the Hessian approximation strategy.}

\subsection{The L-BFGS algorithm}
The Limited-memory BFGS algorithm is a variation of the  Broyden–Fletcher–Goldfarb–Shanno (BFGS) optimization scheme that reduces the amount of memory storage and operations per step of the original algorithm. For the sake of completeness, we briefly review the main points of L-BFGS in the following and refer the reader to \cite{nocedal1999numerical} for more details. L-BFGS belongs to the family of quasi-Newton methods, a class of descent-direction unconstrained optimization schemes that, given an objective function $\varphi$, uses a search direction of the form $\b d_k = -B_k \nabla \varphi(\b x_k)$, with $B_k$ positive definite, to compute the new  approximate minimizer for $\min \varphi$ as $\b x_{k+1}=\b x_k + \alpha_k \b d_k$, with $\alpha_k$ chosen through a suitable line search step. The standard first-order gradient descent method is obtained for $B_k=I$ for all $k$, while the Newton method is obtained by choosing $B_k = \nabla^2\varphi(\b x_k)^{-1}$, the inverse of the Hessian at $\b x_k$. Rather than inverting the Hessian, which can be computationally prohibitive, BFGS computes an approximation of $\nabla^2\varphi(\b x_k)^{-1}$ by performing a rank-2 correction of the previous approximation $B_k = B_{k-1}+R_{k-1}$, with the parameters in the rank-2 matrix $R_{k-1}$ chosen to ensure that (a) $B_k$ is positive definite, and (b) $B_k$ satisfies the secant equation with respect to approximation points $\b x_k$ and $\b x_{k-1}$. This update rule brings down the $\mathcal O(n_F^3)$ cost of Newton's scheme to  $\mathcal O(n_F^2)$ and, more importantly, avoids the computations of second derivatives. To further reduce the cost per step, L-BFGS introduces a ``history parameter'' $m$ and, starting from $B_0=\gamma_0 I$, it updates $B_k$ only for $m$ steps and then resets $B_{k}$ to a multiple of the identity $B_k=\gamma_k I$, every $m$ steps. This operation allows one to further reduce cost and memory storage of the method to $\mathcal O(mn_F)$, which effectively coincides with $\mathcal O(n_F)$ when $m\ll n_F$. In our experiments, we set $m=10$. The pseudocode for L-BFGS is illustrated in Algorithm~\ref{alg:LBFGS}.

In order to apply the L-BFGS approach to the constrained problems \eqref{eq:P-downgrading}, \eqref{eq:P-addition-weighted}, and \eqref{eq:P-tuning}, we modify the objective function  by introducing a logarithmic barrier for the inequality constraints, following a standard Interior-Point method approach (see e.g.\ \cite{byrd1999interior,byrd2000trust}). In \eqref{eq:P-tuning}, for example, the objective function  is modified into
\begin{equation}\label{eq:log-barrier}
  \varphi_\mu(X):= -\varphi_A(X) + \mu \sum_{ij}\big\{\log(U_{ij}-X_{ij})+\log(X_{ij}-A_{ij})\big\}, \qquad \text{with } \mu>0  \, .
\end{equation}
L-BFGS is then applied to the unconstrained problem $\min \varphi_\mu$, and the  parameter $\mu$ is reduced throughout the L-BFGS iterations so that the solution of the approximated problem \eqref{eq:log-barrier} approaches that of \eqref{eq:P-tuning} as the method approaches convergence. In our experiments, the above  Interior-Point method approach with L-BFGS is run by means of Matlab's \texttt{fmincon} function, with optimization parameters \texttt{HessianApproximation}=\texttt{lbfgs} and \texttt{HistorySize}=\texttt{10}.
\begin{algorithm}[t] \small
	\caption{Pseudocode of L-BFGS for unconstrained optimization problem $\min_{\b x} \varphi(\b x)$\label{alg:LBFGS}}
	\begin{algorithmic}[1]
	\Procedure{lbfgs}{$\b x_0$, $\gamma_0$, $m$, $\epsilon$, $\mathrm{maxiter}$}
        \State $B_0 = \gamma_0 I$
        \For{$k=0,1,\dots,\mathrm{maxiter}$}
		\State $\b d_k = -B_k\nabla\varphi(\b x_k)$
        \State $\alpha_k \gets$ line search using $\{\varphi,\b x_k,\b d_k\}$
        \State $\b x_{k+1}=\b x_k + \alpha_k \b d_k$
        \If{$\|\nabla \varphi(\b x_{k+1})\|<\epsilon$}
		 \State \textbf{break}
		\EndIf
        \If{$k+1$ is a multiple of $m$}
		 \State $\gamma_{k+1}\gets$ scalar approximation of $B_k$
        \State $B_{k+1}=\gamma_{k+1}I$
        \Else
        \State $B_{k+1}\gets$ update $B_k$ using the L-BFGS rule
		 \EndIf
		\EndFor
		\State \Return $\b x_{k+1}$
		\EndProcedure
	\end{algorithmic}
\end{algorithm}
\color{black}
\subsection{Gradient approximation via Krylov methods}\label{sec:krylov}
We now look at the gradient of $\varphi_A$, \upd{for a Fr\'echet differentiable $f$}. Let us denote by $\mathrm{ind}:F\rightarrow \{1,\dots,n_F\}$ an ordering map on the set $F$ and observe that the derivative with respect to the $ij$th component of the matrix $X$ is $\partial_{ij} f(A+X)=L_f(A+X, \b{1}_i\b{1}_j^T)$, \upd{where $\b 1_i$ denotes the indicator vector of the node $i$, $(\b 1_i)_j=1$ if $i=j$ and zero otherwise, and}  $L_f(A+X, \b 1_i\b 1_j^T)$ indicates the Fréchet derivative of $f$ at $A+X$, applied to the matrix $\b 1_i\b1_j^T$. Moreover, $L_f(A+X,\b 1_i\b1_j^T+\b 1_j\b1_i^T)= L_f(A+X, \b{1}_i\b{1}_j^T) + L_f(A+X, \b{1}_j\b{1}_i^T)$ and, since $A+X$ is symmetric, it holds $L_f(A+X, \b 1_i\b1_j^T)=L_f(A+X, \b 1_j\b1_i^T)^T $. Putting it all together we have that
\begin{equation}\label{eq:grad}
    \left(\nabla \varphi_A(X)\right)_{\mathrm{ind}(i, j)}=2\, \Tr(L_f(A+X, \b 1_i\b 1_j^T)),\qquad \forall (i,j)\in F.
\end{equation}
\upd{In recent work by Schweitzer \cite{schweitzer23}, it has been shown the following identity
\begin{equation}\label{eq:marcel}
 \Tr(L_f(A+X, \b 1_i\b 1_j^T)) = f'(A^T+X^T)_{ij}= f'(A+X)_{ij},
\end{equation}
where $f'$ denotes the first derivative of $f$. Equation \eqref{eq:marcel} is of key importance because it enables us to simplify the calculation of the gradient from computing $\mathcal O(n_F^2)$ actions of the Fr\'echet derivative to evaluating $\mathcal O(n_F^2)$ entries of a single matrix function. Moreover, if the quantities $f'(A)_{ij}$ are already given, then we can approximate the difference $f'(A+X)_{ij}-f'(A)_{ij}$ with Algorithm~\ref{alg:fun-update}. 
The procedure for evaluating the gradient, for a general $f$, is reported in Algorithm~\ref{alg:gradient}. The cost of the latter is the one of Algorithm~\ref{alg:fun-update} plus extracting $n_F$ entries from the low-rank matrix $\b{\mathcal U}_X\widetilde{\Delta}\b{\mathcal U}_X^*$; under the assumption that matvecs with $A$ cost $\mathcal O(n)$, that $X$ has rank always bounded by $r\leq n_F$, and that Algorithm~\ref{alg:fun-update} takes $\mathrm{it}$ iterations to converge (so that $\widetilde{\Delta}\in\mathbb R^{(r\cdot \mathrm{it})\times (r\cdot \mathrm{it})}$), we get an overall complexity of $\mathcal O((n+n_F)r^2\mathrm{it}^2+r^3\mathrm{it}^4)$.

Note that, when  $f(z)=e^z=f'(z)$, the evaluation of the objective function and of the gradient are based on the same Krylov subspace, i.e., we can compute both by a single execution of the Arnoldi algorithm. More specifically, in the case of the exponential function, we rely on   Algorithm~\ref{alg:fun-update} to both compute the gradient and the objective function $\varphi_A(X)$; the latter requires the quantity $\Tr(\widetilde \Delta )$ that has an additional cost of only $\mathcal O(r \cdot \mathrm{it})$ flops. }
\begin{algorithm} \small
	\caption{Approximation of $\nabla \varphi_A(X)$\label{alg:gradient}}
	\begin{algorithmic}[1]
		\Procedure{gradient\_eval}{$A$, $\b{ U}_X$, $B_X$, $\{f'(A)_{ij}\}_{(i,j)\in\mathrm{ind}^{-1}(\{1,\dots, n_F\})}$, $f'$, $\ell$, $\epsilon$}
		\State $[\b{\mathcal U}_X, \widetilde{\Delta}]\gets \textsc{fun\_update}(A, \b{U}_X, B_{X},f', \ell,\epsilon)$	
		\For{$h=1,\dots, n_F$} 
		\State $(i, j) \gets \mathrm{ind}^{-1}(h)$
		\State $\Delta_h\gets \b{\mathcal U}_X(i,\ :)\cdot \widetilde \Delta \cdot \b{\mathcal U}_X(:,\ j)$		
		\State $\nabla \varphi_A(X)_{h}\gets 2(f'(A)_{ij} + \Delta_{h})$
		\EndFor
		\State \Return $\nabla \varphi_A(X)$
		\EndProcedure
	\end{algorithmic}
\end{algorithm}

\subsection{Hessian evaluation via Krylov methods}
\upd{By taking the partial derivatives of \eqref{eq:marcel} we get the following expression for the Hessian's entries:
\begin{equation}\label{eq:hessian}
	\left(H\varphi_A(X)\right)_{\mathrm{ind}(i,j), \mathrm{ind}(h,k)}= 2\left(L_{f'}(A+X, \b 1_i\b 1_j^T)\right)_{hk}\qquad \forall (i,j),(h,k)\in F.
\end{equation}
In particular, \eqref{eq:hessian} tells us that computing the Hessian requires extracting $\mathcal O(n_F)$ entries from $\mathcal O(n_F)$ Fréchet derivatives along rank $1$ directions. Fortunately, the rank $1$ property of the direction implies the low-rank approximability of $L_{f'}(A+X, \b 1_i\b 1_j^T)$ that in turn enables us to leverage an efficient Krylov subspace technique \cite{kressner2019bivariate}, as discussed next.} 

\upd {To simplify the exposition we temporarily replace $f'$ with $f$ and we describe how to efficiently evaluate quantities  of the form $L_f(M, \b{1}_i\b{1}_j^T)$, for a given symmetric  matrix $M$ and a given function $f$.} 
The evaluation of the Fr\'echet derivative in a certain direction can be recast as evaluating the function of a specific augmented matrix. More precisely, applying the well-known formula in \cite[Theorem 2.1]{Roy96} to our framework, yields
\begin{equation}\label{eq:frechet:fA}
f\left( \begin{bmatrix}M & \b{1}_i\b{1}_j^T\\ 0 & M\end{bmatrix} \right) = \begin{bmatrix} f(M) & L_f(M,\b{1}_i\b{1}_j^T) \\ 0 & f(M) \end{bmatrix}
\end{equation}
so that we can look at extracting the $(1,2)$ sub-block of \eqref{eq:frechet:fA}.
Since  $$\begin{bmatrix} 0 & L_f(M,\b{1}_i\b{1}_j^T) \\ 0 & 0 \end{bmatrix}=f\left( \begin{bmatrix}M & 0\\ 0 & M\end{bmatrix} + \begin{bmatrix}0 & \b{1}_i\b{1}_j^T\\ 0 & 0\end{bmatrix}\right)-f\left(\begin{bmatrix}M & 0\\ 0 & M\end{bmatrix}\right)$$
and $\begin{bmatrix}0 & \b{1}_i\b{1}_j^T\\ 0 & 0\end{bmatrix}$ is of rank $1$, we expect $L_f(M,\b{1}_i\b{1}_j^T)$ to be well approximated by a low-rank matrix.  This property is exploited in  \cite[Algorithm 2]{kressner2019bivariate}, where a projection method that makes use of tensorized Krylov subspaces has been proposed. The latter incrementally builds   orthonormal bases $\b{\mathcal U}_m,\b{\mathcal V}_m$ for $\mathcal K_m(M, \b 1_i)$ and  $\mathcal K_m(M, \b 1_j)$, respectively, by means of two Arnoldi processes. The  associated Arnoldi relations
$$
 M \b{\mathcal U}_m = \b{\mathcal U}_m \mathcal H_m + \b U_{m+1} H_{m+1,m} \b 1_m^T,\qquad M \b{\mathcal V}_m = \b{\mathcal V}_m \mathcal G_m + \b U_{m+1} G_{m+1,m} \b 1_m^T, 
$$
directly provide the expression of the projected augmented matrix 
$$
\begin{bmatrix}  \b{\mathcal U}_m^*& 0 \\ 0 & \b{\mathcal V}_m^* \end{bmatrix}\begin{bmatrix} M &\b{1}_i\b{1}_j^T \\ 0 & M \end{bmatrix}\begin{bmatrix}  \b{\mathcal U}_m& 0 \\ 0 & \b{\mathcal V}_m \end{bmatrix}=\begin{bmatrix}  \mathcal H_m& \b{1}_i\b{1}_j^T \\ 0 & \mathcal G_m \end{bmatrix}.
$$
Thus, the method computes the quantities
\begin{equation}\label{eq:frecht-approx}
    L_{f,m}^{(i,j)}:=\b{\mathcal U}_m\widetilde{L}_{f,m}^{(i,j)}\b{\mathcal V}_m^*,\qquad \widetilde{L}_{f,m}^{(i,j)}:=f\left(\begin{bmatrix}  \mathcal H_m& \b{1}_i\b{1}_j^T \\ 0 & \mathcal G_m \end{bmatrix}\right)_{(1,2)},
\end{equation}
where the subscript $(1,2)$ refers to the extraction of the $(1,2)$ sub-block, as an approximation of $L_f(M,\b{1}_i\b{1}_j^T)$. 
The method then stops when the heuristic stopping criterion 
$$
\norm{\widetilde{L}_{f,m}^{(i,j)}-\begin{bmatrix}\widetilde{L}_{f,m-\ell}^{(i,j)}& 0 \\ 0 &0\end{bmatrix}}_2\leq \epsilon
$$  
is verified, for a prescribed tolerance $\epsilon$ and a positive integer $\ell$. In our implementation we set $\ell=2$. 
\upd{For an alternative and more reliable stopping criterion see \cite[Section 5]{kandolf2021computing}.}
The full procedure is reported in Algorithm~\ref{alg:frechet-update}.
\begin{algorithm}[t] \small
	\caption{Approximation of $L_f(M,\b 1_i\b 1_j^T)$\label{alg:frechet-update}}
	\begin{algorithmic}[1]
	\Procedure{frechet\_eval}{$M$, $i$, $j$, $f$, $\ell$, $\epsilon$}
		\For{$m = 1, \ldots, m_{\max}$}
		\State \begin{minipage}{0.8\textwidth} Compute (incrementally) the Arnoldi relation for $\mathcal K_{m}(M,\b 1_i),\mathcal K_{m}(M,\b 1_j)$ \upd{by means of the Arnoldi method}; store 
		 $\b{\mathcal U}_{m},\b{\mathcal V}_{m},\mathcal H_m$ and $\mathcal G_m$
		 \end{minipage}
		 \State $\widetilde{L}_{f,m}^{(i,j)}\gets f\left(\begin{bmatrix}  \mathcal H_m& \b{1}_i\b{1}_j^T \\ 0 & \mathcal G_m \end{bmatrix}\right)_{(1,2)}$
		 \If{$m> \ell$ and $\norm{\widetilde{L}_{f,m}^{(i,j)}-\begin{bmatrix}\widetilde{L}_{f,m-\ell}^{(i,j)}& 0 \\ 0 &0\end{bmatrix}}_2\leq \epsilon$}
		 \State \textbf{break}
		 \EndIf
		\EndFor
		\State \Return $\b{\mathcal U}_{m},\b{\mathcal V}_{m},\widetilde{L}_{f,m}^{(i,j)}$
		\EndProcedure
	\end{algorithmic}
\end{algorithm}

We point out that, the approximation error associated with the sequence $L_{f,m}^{(i,j)}$, $m=1,2,\dots$, decays at least as the best polynomial approximation error of $f'$ on the convex hull  of the spectrum of $M$, which we denote by $\Pi$. More precisely, a direct consequence of \cite[Corollary 1]{kressner2019bivariate} is the following bound:
\[
\norm{L_f(M, \b{1}_i\b{1}_j^T)-L_{f,m}^{(i,j)} )}_F\leq 2\min_{p\in\mathcal P_{m-1}}\max_{z\in \Pi} |f'(z)-p(z)|.
\]
\upd{Note that, the cost analysis of Algorithm~\ref{alg:frechet-update} is very similar to the one of Algorithm~\ref{alg:fun-update}. In particular, under the assumptions that matvecs with $M$ cost $\mathcal O(n)$, and that the Arnoldi procedure takes $\mathrm{it}$ iterations before detecting convergence, Algorithm~\ref{alg:frechet-update} costs $\mathcal O(n\cdot \mathrm{it}^2+\mathrm{it}^4)$.}

\subsubsection{Multiple evaluations of $L_f(M,\b 1_i\b 1_j^T)$}
Evaluating \eqref{eq:hessian} requires to approximate the quantities $L_{f'}(A,\b 1_i\b 1_j^T)$ for all edges $(i,j)\in F$, with $|F|=n_F\ll n$.  In principle, running Algorithm~\ref{alg:frechet-update} $n_F$ times (on each pair $(i,j)\in F$) performs the sought evaluation. On the other hand, it is possible to  enhance the efficiency by avoiding redundant computations due to the repetition of the same nodes in edges of $F$ and thus the same Krylov subspaces.  Denote by $V(F)$ the set of nodes that are linked by the edges in $F$, i.e., $V(F):=\{i\in V: \exists j\in V \text{ such that } (i,j)\in F\}$ and for any such node $i\in V(F)$ let $V_i(F)$ be the set of nodes  that are connected to $i$ via an edge in $F$, i.e., $V_i(F):=\{j\in V: (i,j)\in F\}\subseteq V(F)$. We proceed as follows:
\begin{enumerate}[(i)]
\item For each $i\in V(F)$ we compute and store the Arnoldi relation 
\[M \b{\mathcal U}_{m_i}^{(i)} = \b{\mathcal U}_{m_i}^{(i)} \mathcal H_{m_i}^{(i)} + \b{U}_{m_i+1}^{(i)} H_{m_i+1,m_i}^{(i)} \b 1_{m_i}^T
\]
for $\mathcal K_{m_i}(M, \b 1_i)$ where $m_i$ is such that $\norm{\widetilde{L}_{f,m_i}^{(i,j)}-\left[\begin{smallmatrix}\widetilde{L}_{f,m_i-\ell}^{(i,j)}& 0 \\ 0 &0\end{smallmatrix}\right]}_2\leq \epsilon$, for all $j\in V_i(F)$. 
    
\item While doing $(i)$, for each pair $(i,j)\in F$, we store $\widetilde{L}_{f,m_{(i,j)}}^{(i,j)}$ where  $m_{(i,j)}$ is the smallest integer such that $\norm{\widetilde{L}_{f,m_{(i,j)}}^{(i,j)}-\left[\begin{smallmatrix}\widetilde{L}_{f,m_{(i,j)}-\ell}^{(i,j)}& 0 \\ 0 &0\end{smallmatrix}\right]}_2\leq \epsilon$.  Note that, $m_{(i,j)}\leq m_i$ which may yield a cheaper trace evaluation for that particular $(i,j)\in F$.
\end{enumerate}
The procedure which implements these enhancements is reported in Algorithm~\ref{alg:multiple-frechet-eval}. 

\upd{Let us denote by $n_{V}:= |V(F)|$ and $\mathrm{it}=\max_{(i,j)\in F}m_{(i,j)}$ and assume that matvecs with $M$ cost $\mathcal O(n)$. Then, the complexity of Algorithm~\ref{alg:multiple-frechet-eval} is determined by $n_V$ times the one of Algorithm~\ref{alg:frechet-update}, i.e., $\mathcal O(n_V(n\cdot \mathrm{it}^2+ \mathrm{it}^4))$. We remark that Algorithm~\ref{alg:multiple-frechet-eval} requires to store $\mathcal O(n_V\cdot \mathrm{it})$ vectors of length $n$ to represent all the Krylov bases; this
	might not be feasible for a large value of $n_V$, i.e., a large search space $F$.}  

 \upd{Finally, the procedure that evaluates the Hessian of $\varphi_A(X)$ is reported in Algorithm~\ref{alg:hessian}.
The latter consists in one call to Algorithm~\ref{alg:multiple-frechet-eval} and extracting $ n_F$ entries from a matrix of rank $r\cdot \mathrm{it}$, $\mathcal O(n_F)$ times; this yields a complexity estimate of $\mathcal O(n_V(n\cdot \mathrm{it}^2+ \mathrm{it}^4) + n_F^2\mathrm{it}^2)$.
}
\begin{algorithm}[t] \small
	\caption{Approximation of $L_f(M,\b 1_i\b 1_j^T)$ $\forall (i,j)\in F$\label{alg:multiple-frechet-eval}}
	\begin{algorithmic}[1]
	\Procedure{multiple\_frechet\_eval}{$M$, $F$, $f$, $\ell$, $\epsilon$}
	\State NC $\gets F$, $\qquad m=1$ \Comment{NC is the set of not converged edges}
		\While{NC$\neq\emptyset$}
		\For{$i\in V(F)$}
		\If{$\exists j\in V_i(F)$ such that $(i,j)\in$ NC}
		\State \begin{minipage}{0.8\textwidth} Compute (incrementally) the Arnoldi relation for $\mathcal K_{m}(M,\b 1_i)$ \upd{by means of the Arnoldi method}; store 
		 $\b{\mathcal U}_{m}^{(i)}$ and $\mathcal H_m^{(i)}$ 
		 \end{minipage}
		 \EndIf
		 \EndFor
		 \For{$(i,j)\in$ NC}
		 \State $\widetilde{L}_{f,m}^{(i,j)}\gets f\left(\begin{bmatrix}  \mathcal H_m^{(i)}& \b{1}_i\b{1}_j^T \\ 0 & \mathcal H_m^{(j)} \end{bmatrix}\right)_{(1,2)}$
		 \If{$m> \ell$ \textbf{and} $\norm{\widetilde{L}_{f,m}^{(i,j)}-\left[\begin{smallmatrix}\widetilde{L}_{f,m-\ell}^{(i,j)}& 0 \\ 0 &0\end{smallmatrix}\right]}_2\leq \epsilon$}
		 \State NC $\gets$ NC$\setminus (i,j)$, $\qquad m_{(i,j)}\gets m$
		 \EndIf
		 \EndFor
		 \State $m\gets m+1$
		 \EndWhile
		\State \Return $\b{\mathcal U}_{m}^{(h)}$ for all nodes $h\in V(F)$ and $\widetilde L_{f,m_{(i,j)}}^{(i,j)}$  for all edges $(i,j)\in F$
		\EndProcedure
	\end{algorithmic}
\end{algorithm}

\begin{algorithm} \small
	\caption{Approximation of $H \varphi_A(X)$\label{alg:hessian}}
	\begin{algorithmic}[1]
		\Procedure{hessian\_eval}{$A$, $\b{ U}_X$, $B_X$, $F$, $f'$, $\ell$, $\epsilon$}
		\State $\{\b{\mathcal U}_{\mathrm{ind}(i,j)}, \widetilde L^{(i,j)}\}_{(i,j)\in F}\gets\textsc{multiple\_frechet\_update}(A+\b{ U}_XB_X\b{U}^*, F, f',\ell)$
		\For{$s=1,\dots, n_F$} 
				\State $(i,\ j) \gets \mathrm{ind}^{-1}(s)$
		\For{$t=s,\dots, n_F$}
		\State $(h,\ k) \gets \mathrm{ind}^{-1}(t)$		
		\State $H \varphi_A(X)_{st}\gets 2\  \b{\mathcal U}_s(h,\ 1:m_{(i,j)})\  \widetilde L^{(i,j)}\ \b{\mathcal U}_s(k,\ 1:m_{(i,j)})^*$
		\State $H \varphi_A(X)_{ts}\gets H \varphi_A(X)_{st}$
		\EndFor
		\EndFor
		\State \Return $H \varphi_A(X)$
		\EndProcedure
	\end{algorithmic}
\end{algorithm}


\section{Numerical experiments with unweighted graphs}\label{sec:experiments}
We test the performance of \textsc{greedy\_krylov\_break} and \textsc{greedy\_krylov\_make}, \upd{introduced in Section~\ref{sec:competitors}}, with respect to their effectiveness in manipulating the graph natural connectivity, i.e. $f(z)=e^z$, and their running time on 22 real-world unweighted networks. Details about the networks' size are reported in Table~\ref{fig:datasets}; \upd{in particular, the corresponding adjacency matrices are of size $|V|\times |V|$ and have at most $2|E|$ nonzero entries.} Those listed on the left-hand side of Table~\ref{fig:datasets}  include social networks of geolocated reciprocated Twitter mentions within UK cities (Cardiff, Edinburgh), coauthorship networks (ca-AstroPh, ca-CondMat, ca-HephTh, netscience), a protein-protein interactions network (yeast) and a public transports network (London). All these networks are publicly available via public repositories, as reported in \cite{cipolla2021nonlocal,grindrod2016comparison,Pajek,snap}.  All the  networks listed on the right-hand side of Table~\ref{fig:datasets}  are road networks of different cities in the world \cite{roadnetworks}. Our implementation is written using MATLAB and is available at the public repository \url{https://github.com/COMPiLELab/krylov_robustness}, together with all the datasets above.


In the proposed experiments we compare with  state-of-the-art methods \textsc{miobi} and \textsc{eigenv},  that have been recalled in Section~\ref{sec:competitors}. In particular, \textsc{miobi} uses $25$ eigenpairs to compute the approximate trace variation as described in Section~\ref{sec:miobi} and the search spaces \ref{strat:exhaustive}, \ref{strat:degree} for problems \ref{eq:P-removal} and \ref{eq:P-addition-binary}, respectively. \upd{If not stated otherwise, the $f$-connectivity is considered with respect to be  the matrix exponential function, i.e., $f=\exp$.}

To assess the impact of the various methods on the natural connectivity of a network we consider the \emph{magnitude of the relative trace variation} that, given the returned modification of the adjacency matrix $X$, we define as:
$$\Delta T (X):=\frac{|\Tr(f(A+X))-\Tr(f(A))|}{|\Tr(f(A))|}.$$
\upd{To obtain an estimate of the denominator $\Tr(f(A))$ we have employed the  stochastic trace estimator  \textsc{hutch++}~\cite{hutch++} combined with the \texttt{expmv} algorithm from~\cite{almohy} to evaluate the action of the matrix exponential on Rademacher random vectors.}

Finally, to evaluate the scalability of the approaches we report their computational times in seconds. \upd{The latter do not include the time spent for estimating $\Tr(f(A))$ at the beginning, as this operation is not required by the greedy procedures.}

\upd{The experiments have been performed on a laptop with a dual-core Intel
Core i7-7500U 2.70 GHz CPU, 256KB of level 2 cache, 16 GB of RAM, and operating system Ubuntu 22.04.2. The algorithms are
implemented in MATLAB and tested under MATLAB2022b, with MKL BLAS version 2019.0.3
utilizing both cores.} 

\begin{table}[t]
	\caption{Number of vertices and edges of the unweighted graphs used for the numerical tests. On the left: social, collaboration, transportation, and PPI networks; on the right: graphs representing road networks.}\label{fig:datasets}
\vspace{.2cm}
	
	\centering
	\begin{tabular}{lcc}
	\toprule	
	\textbf{Dataset}& $|V|$&$|E|$\\
	\midrule
	Cardiff&$2685$&$4444$\\
	CollegeMsg&$1893$&$13835$\\
	Edinburgh&$1645$&$2146$\\
	as\_735&$6474$&$12572$\\
	ca-AstroPh&$17903$&$19972$\\
	ca-CondMat&$21363$&$91286$\\
	ca-HephTh&$8638$&$24806$\\
	London&$369$&$430$\\
	netscience&$379$&$914$\\
	socEpinions1&$75877$&$405739$\\
	yeast&$2224$&$6609$\\
	\bottomrule
\end{tabular}
\qquad\qquad 
	\begin{tabular}{lcc}
		\toprule	
		\textbf{Dataset}& $|V|$&$|E|$\\
		\midrule
		Anaheim&$416$&$634$\\
		Austin&$7388$&$10591$\\
		Barcelona&$930$&$1798$\\
		Birmingham&$14578$&$20913$\\
		ChicagoRegional&$12979$&$20627$\\
		DC&$9522$&$14807$\\
		Hawaii&$21774$&$26007$\\
		Philadelphia&$13389$&$21246$\\
		RhodeIsland&$51642$&$66650$\\
		Rome&$3353$&$4831$\\
		Sydney&$32956$&$38787$\\
		\bottomrule
	\end{tabular}
\end{table}
\subsection{Downgrading for unweighted graphs}
As a first experiment, we measure the quantity $\Delta T$ when solving problem \eqref{eq:P-removal} with a fixed budget of $k=50$ edges \upd{to be removed}. The parameter $q$, used by the method \textsc{greedy\_krylov\_break} to determine its search space, is set to the value $250$. 

The performances of \textsc{greedy\_krylov\_break}, \textsc{miobi}, and \textsc{eigenv} are compared over both  road and general networks. The results reported in the left part of Table~\ref{tab:break_road_misc} show that \textsc{miobi} and \textsc{greedy\_krylov\_break} always outperform \textsc{eigenv} on road networks and our \textsc{greedy\_krylov\_break} achieves the best score on $6$ out of $11$ case studies.
Also, for general graphs, \textsc{miobi} and \textsc{greedy\_krylov\_break} provide the best scores although the results reported in the right part of Table~\ref{tab:break_road_misc} show a balanced situation: on $7$ out of $11$ case studies, the difference between the scores of the methods is less than 2\%. The most evident gain of the top method is measured for the medium-size graph \texttt{ca-HephTh} and the small graph \texttt{netscience}. In view of the significantly lower costs of \textsc{miobi} and \textsc{eigenv} (see Section~\ref{sec:exp-budget}), these results suggest that \textsc{greedy\_krylov\_break} can be a valid competitor for the road networks dataset only.

\subsubsection{Trace reduction and scalability with respect to the budget size}\label{sec:exp-budget}
Now we consider a second numerical test where we let the budget size $k$ range in the set $\{10\cdot j\}$, $j=1,\dots, 10$, and we measure both the relative trace variation and the time consumption of the methods. Further, we  investigate how  the parameter $q$, \upd{that determines the size of the search space}, affects the performance of \textsc{greedy\_krylov\_break} by considering three implementations of this method for $q=50,250,\min\{1000,|E|-k\}$. As case studies, we select $6$ road networks: \texttt{Anaheim}, \texttt{Birmingham}, \texttt{ChicagoRegional}, \texttt{Hawaii}, \texttt{RhodeIsland}, and \texttt{Rome}. Figure~\ref{fig:break_score} reports the magnitude of the relative trace variations attained by the five methods, as the budget increases. The method \textsc{greedy\_krylov\_break} with the largest search space attains the highest scores on all the examples apart from \texttt{Birmingham}, where the returned trace variation is comparable with the one of \textsc{miobi}. There is no clear winner between \textsc{miobi} and \textsc{greedy\_krylov\_break} with $q=500$, while  \textsc{greedy\_krylov\_break} with $q=50$ and \textsc{eigenv} always provide the $4$th and the $5$th scores.

\begin{figure}
	\centering
	\begin{tikzpicture}[scale=0.85]
	\begin{groupplot}[group style={group size= 3 by 2, vertical sep=1.5cm}, legend style={					
			font=\scriptsize,
			at={(3.69,-1.75)},
			legend columns = -1,
			cells={anchor=west}
		}, width=.4\textwidth]
		\nextgroupplot[title=Anaheim, scaled y ticks=base 10:1, ymode=log]
		\addplot table[x index = 0, y index = 1] {break_anaheim_time.dat};\addlegendentry{\textsc{greedy\_krylov\_break}($q=50$)};
		\addplot table[x index = 0, y index = 2] {break_anaheim_time.dat};\addlegendentry{\textsc{greedy\_krylov\_break}($q=250$)};
		\addplot table[x index = 0, y index = 3] {break_anaheim_time.dat};\addlegendentry{\textsc{greedy\_krylov\_break}($q=1000$)};
		\addplot table[x index = 0, y index = 4] {break_anaheim_time.dat};\addlegendentry{\textsc{miobi}};
		\addplot[mark = triangle, green] table[x index = 0, y index = 5] {break_anaheim_time.dat};\addlegendentry{\textsc{eigenv}};
		\nextgroupplot[title=Birmingham, ymode=log]
		\addplot table[x index = 0, y index = 1] {break_birmingham_time.dat};
		\addplot table[x index = 0, y index = 2] {break_birmingham_time.dat};
		\addplot table[x index = 0, y index = 3] {break_birmingham_time.dat};
		\addplot table[x index = 0, y index = 4] {break_birmingham_time.dat};
		\addplot[mark = triangle, green] table[x index = 0, y index = 5] {break_birmingham_time.dat};
		\nextgroupplot[title=ChicagoRegional, ymode=log]
		\addplot table[x index = 0, y index = 1] {break_chicago_time.dat};
		\addplot table[x index = 0, y index = 2] {break_chicago_time.dat};
		\addplot table[x index = 0, y index = 3] {break_chicago_time.dat};
		\addplot table[x index = 0, y index = 4] {break_chicago_time.dat};
		\addplot[mark = triangle, green] table[x index = 0, y index = 5] {break_chicago_time.dat};
		\nextgroupplot[title=Hawaii, ylabel= Computational time \upd{(sec.)},	every axis y label/.append style={at=(ticklabel cs:1.1)}, xlabel = Budget of edges to be removed ($k$), every axis x label/.append style={at=(ticklabel cs:1.75)}, ymode=log]
		\addplot table[x index = 0, y index = 1] {break_hawaii_time.dat};
		\addplot table[x index = 0, y index = 2] {break_hawaii_time.dat};
		\addplot table[x index = 0, y index = 3] {break_hawaii_time.dat};
		\addplot table[x index = 0, y index = 4] {break_hawaii_time.dat};
		\addplot[mark = triangle, green] table[x index = 0, y index = 5] {break_hawaii_time.dat};
		\nextgroupplot[title=RhodeIsland, ymode=log]
		\addplot table[x index = 0, y index = 1] {break_rhodeisland_time.dat};
		\addplot table[x index = 0, y index = 2] {break_rhodeisland_time.dat};
		\addplot table[x index = 0, y index = 3] {break_rhodeisland_time.dat};
		\addplot table[x index = 0, y index = 4] {break_rhodeisland_time.dat};
		\addplot[mark = triangle, green] table[x index = 0, y index = 5] {break_rhodeisland_time.dat};
		
		\nextgroupplot[title=Rome, ymode=log]
		\addplot table[x index = 0, y index = 1] {break_rome_time.dat};
		\addplot table[x index = 0, y index = 2] {break_rome_time.dat};
		\addplot table[x index = 0, y index = 3] {break_rome_time.dat};
		\addplot table[x index = 0, y index = 4] {break_rome_time.dat};
		\addplot[mark = triangle, green] table[x index = 0, y index = 5] {break_rome_time.dat};
	\end{groupplot}
\end{tikzpicture}
	\caption{Computational times \upd{(seconds)} of the methods for downgrading as the budget increases.}\label{fig:break_time}
\end{figure}
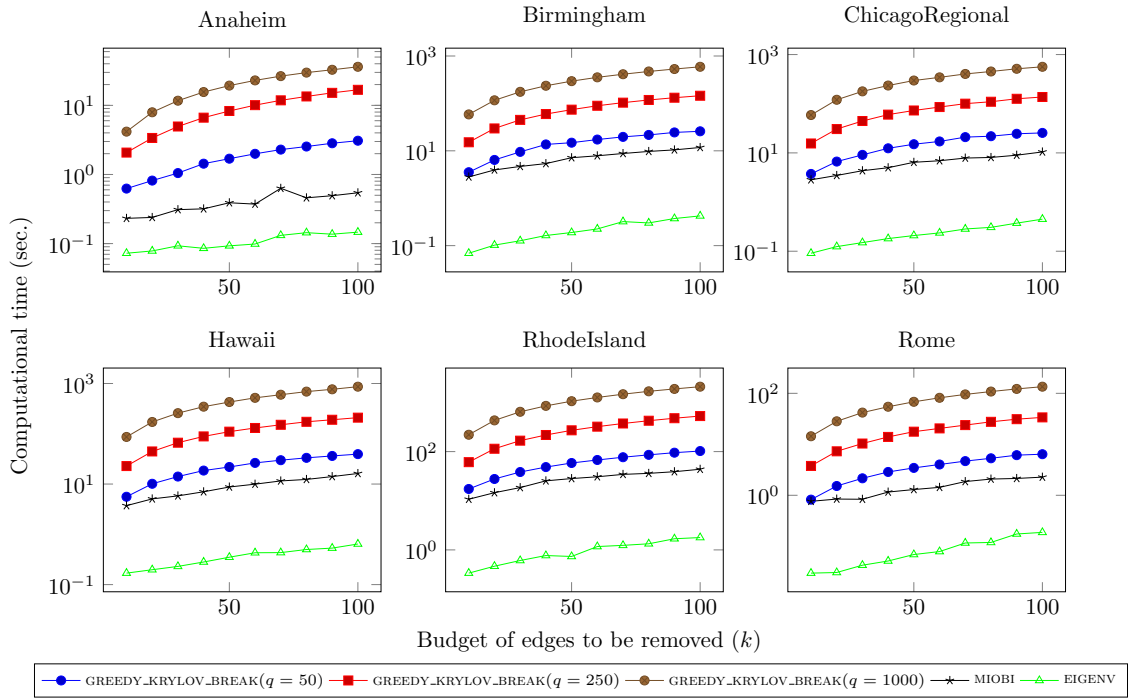

The computational times shown in Figure~\ref{fig:break_time} confirm that the cost of all algorithms has a linear scaling with respect to the parameter $k$. Also, their dependence on $n$ is linear, but the hidden constant determines significantly different running times. In particular, in all case studies the three implementations of \textsc{greedy\_krylov\_break} are the most expensive, then we have \textsc{miobi} and, finally, \text{eigenv} that is the cheapest method. As expected, reducing the parameter $q$ improves the timings of \textsc{greedy\_krylov\_break}, however, in view of the scores in Figure~\ref{fig:break_score}, the convenience of a smaller search space is questionable. Overall, these results suggest that \textsc{greedy\_krylov\_break} is preferable in a scenario where the robustness reduction matters more than the computing time.
\begin{table}[t]
    \centering
    \textbf{Downgrading}
    \resizebox{\textwidth}{!}{\begin{tabular}{lc@{\hskip .5em}c@{\hskip .5em}c@{\hskip .5em}c c@{\hskip 2em} lc@{\hskip .5em}c@{\hskip .5em}c@{\hskip .5em}c}
    \toprule
    & GKB&MIOBI&EIGENV&$\cap$  &&&GKB&MIOBI&EIGENV&$\cap$\\
    \cmidrule{1-5}  \cmidrule{7-11}
    Anaheim&\textbf{0.123}&0.0956&0.0775&6&&Cardiff&0.974&\textbf{0.974}&0.973&43\\
    Austin&0.00863&\textbf{0.00943}&0.00564&6&&CollegeMsg&\textbf{0.773}&0.771&0.771&45\\
    Barcelona&0.0871&\textbf{0.0900}&0.0634&10&&Edinburgh&0.326&\textbf{0.335}&0.240&19\\
    Birmingham&0.00364&\textbf{0.00478}&0.00234&4&&as735&0.965&\textbf{0.966}&\textbf{0.966}&45\\
    ChicagoRegional&\textbf{0.00530}&0.00501&0.00317&4&&AstroPh&\textbf{0.751}&0.751&0.751&49\\
    DC&\textbf{0.00682}&0.00643&0.00417&5&&CondMat&\textbf{0.858}&0.854&0.854&45\\
    Hawaii&0.00273&\textbf{0.00287}&0.00198&7&&HepTh&\textbf{0.958}&0.847&0.847&5\\
    Philadelphia&\textbf{0.00348}&0.00340&0.00236&2&&London&\textbf{0.158}&0.151&0.119&12\\
    RhodeIsland&\textbf{0.00125}&0.00124&0.000752&2&&netscience&0.704&\textbf{0.814}&0.744&18\\
    Rome&\textbf{0.0161}&0.0158&0.0101&3&&Epinions1&0.581&\textbf{0.587}&\textbf{0.587}&41\\
    Sydney&0.00148&\textbf{0.00250}&0.00109&2&&yeast&\textbf{0.878}&0.871&0.865&36\\
    \bottomrule
    \end{tabular}}
    \caption{Magnitude of the relative trace variation obtained with the three methods \textsc{greedy\_krylov\_break (gkb)}, \textsc{miobi}, \textsc{eigenv} considered for the downgrading of unweighted graphs on road networks (left) and general networks (right), with a budget of $k=50$ edges. The column denoted with $\cap$ shows the number of edges that have been commonly chosen  by all the methods.}
    \label{tab:break_road_misc}
\end{table}

	\begin{figure}
	\centering
	\begin{tikzpicture}[scale=.85]
		\begin{groupplot}[group style={group size= 3 by 2, vertical sep=1.5cm}, legend style={					
				font=\scriptsize,
				at={(3.69,-1.75)},
				legend columns = -1,
				cells={anchor=west}
			}, width=.4\textwidth]
				\nextgroupplot[title=Anaheim, scaled y ticks=base 10:1]
		\addplot table[x index = 0, y index = 1] {break_anaheim.dat};\addlegendentry{\textsc{greedy\_krylov\_break}($q=50$)};
		\addplot table[x index = 0, y index = 2] {break_anaheim.dat};\addlegendentry{\textsc{greedy\_krylov\_break}($q=250$)};
		\addplot table[x index = 0, y index = 3] {break_anaheim.dat};\addlegendentry{\textsc{greedy\_krylov\_break}($q=1000$)};
		\addplot table[x index = 0, y index = 4] {break_anaheim.dat};\addlegendentry{\textsc{miobi}};
		\addplot[mark = triangle, green] table[x index = 0, y index = 5] {break_anaheim.dat};\addlegendentry{\textsc{eigenv}};
			\nextgroupplot[title=Birmingham]
			\addplot table[x index = 0, y index = 1] {break_birmingham.dat};
			\addplot table[x index = 0, y index = 2] {break_birmingham.dat};
			\addplot table[x index = 0, y index = 3] {break_birmingham.dat};
			\addplot table[x index = 0, y index = 4] {break_birmingham.dat};
			\addplot[mark = triangle, green] table[x index = 0, y index = 5] {break_birmingham.dat};
				\nextgroupplot[title=ChicagoRegional]
			\addplot table[x index = 0, y index = 1] {break_chicago.dat};
			\addplot table[x index = 0, y index = 2] {break_chicago.dat};
			\addplot table[x index = 0, y index = 3] {break_chicago.dat};
			\addplot table[x index = 0, y index = 4] {break_chicago.dat};
			\addplot[mark = triangle, green] table[x index = 0, y index = 5] {break_chicago.dat};
				\nextgroupplot[title=Hawaii, ylabel= Magnitude of the relative trace variation,	every axis y label/.append style={at=(ticklabel cs:1.1)}, xlabel = Budget of edges to be removed ($k$), every axis x label/.append style={at=(ticklabel cs:1.75)}]
			\addplot table[x index = 0, y index = 1] {break_hawaii.dat};
			\addplot table[x index = 0, y index = 2] {break_hawaii.dat};
			\addplot table[x index = 0, y index = 3] {break_hawaii.dat};
			\addplot table[x index = 0, y index = 4] {break_hawaii.dat};
			\addplot[mark = triangle, green] table[x index = 0, y index = 5] {break_hawaii.dat};
	\nextgroupplot[title=RhodeIsland]
		\addplot table[x index = 0, y index = 1] {break_rhodeisland.dat};
	\addplot table[x index = 0, y index = 2] {break_rhodeisland.dat};
	\addplot table[x index = 0, y index = 3] {break_rhodeisland.dat};
	\addplot table[x index = 0, y index = 4] {break_rhodeisland.dat};
	\addplot[mark = triangle, green] table[x index = 0, y index = 5] {break_rhodeisland.dat};

		\nextgroupplot[title=Rome]
	\addplot table[x index = 0, y index = 1] {break_rome.dat};
	\addplot table[x index = 0, y index = 2] {break_rome.dat};
	\addplot table[x index = 0, y index = 3] {break_rome.dat};
	\addplot table[x index = 0, y index = 4] {break_rome.dat};
	\addplot[mark = triangle, green] table[x index = 0, y index = 5] {break_rome.dat};
   \end{groupplot}
	\end{tikzpicture}
	\caption{Magnitude of the relative trace variation for downgrading as the budget increases.}\label{fig:break_score}
\end{figure}
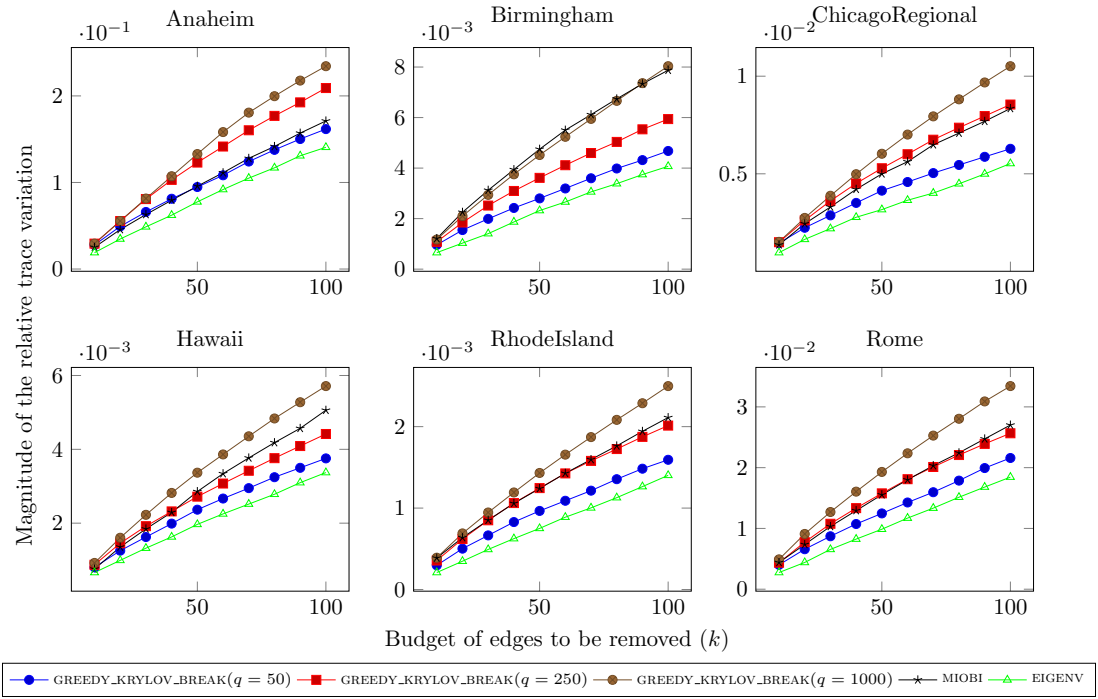

\subsection{Addition for unweighted graphs}
Here we consider analogous tests to those performed in the previous section, for the optimization problem \eqref{eq:P-addition-binary}. This time, we compare \textsc{miobi}  and \textsc{eigenv} with the performance of our \textsc{greedy\_krylov\_make} with $q=\min\{1000, |E|\}$.  In  the left and right parts of Table~\ref{tab:make_road_misc} it is reported the magnitude of the relative trace variation, obtained with \upd{a budget} $k=50$, for road and general networks, respectively. For all road networks, \textsc{greedy\_krylov\_make} is the clear-cut winner and outperforms the second-highest score of a factor  between $1.5$ and $5$. For general networks, \textsc{greedy\_krylov\_make} obtains either the best or near-best score on $10$ out of $11$ examples, although the gain with respect to the competitors is often more limited than for road networks.

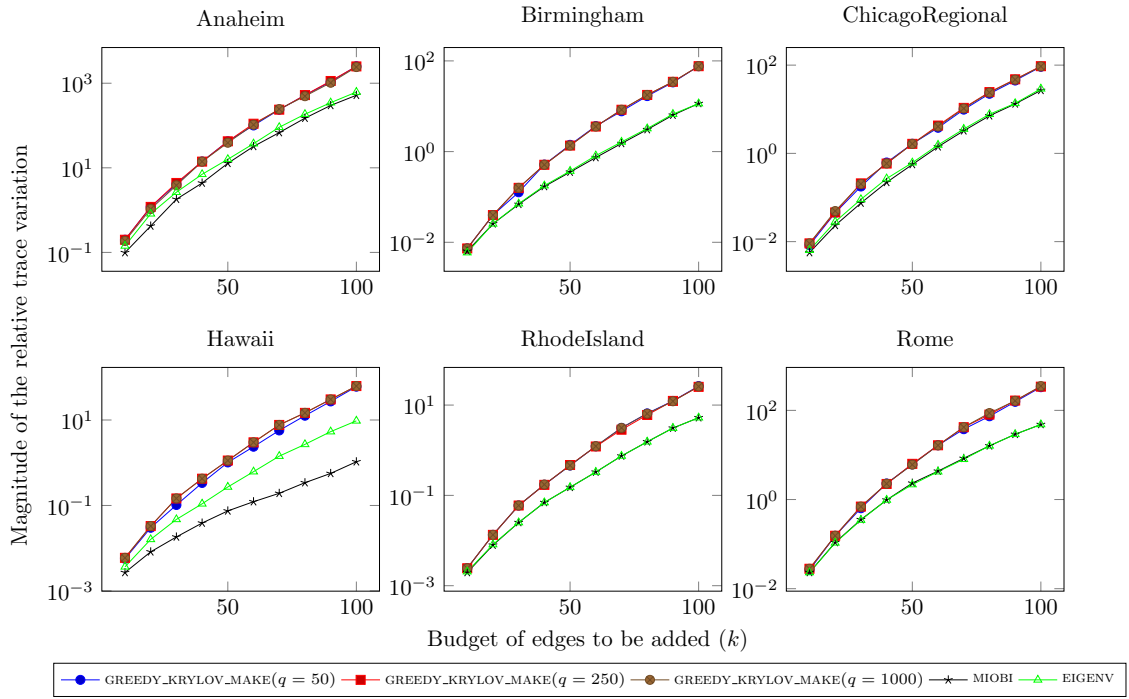
\begin{figure}
	\centering
	\begin{tikzpicture}[scale=.85] 
		\begin{groupplot}[group style={group size= 3 by 2, vertical sep=1.5cm}, legend style={					
				font=\scriptsize,
				at={(3.69,-1.75)},
				legend columns = -1,
				cells={anchor=west}
			}, width=.4\textwidth]
			\nextgroupplot[title=Anaheim, scaled y ticks=base 10:1, ymode = log]
			\addplot table[x index = 0, y index = 1] {make_anaheim.dat};\addlegendentry{\textsc{greedy\_krylov\_make}($q=50$)};
			\addplot table[x index = 0, y index = 2] {make_anaheim.dat};\addlegendentry{\textsc{greedy\_krylov\_make}($q=250$)};
			\addplot table[x index = 0, y index = 3] {make_anaheim.dat};\addlegendentry{\textsc{greedy\_krylov\_make}($q=1000$)};
			\addplot table[x index = 0, y index = 4] {make_anaheim.dat};\addlegendentry{\textsc{miobi}};
			\addplot[mark = triangle, green] table[x index = 0, y index = 5] {make_anaheim.dat};\addlegendentry{\textsc{eigenv}};
			\nextgroupplot[title=Birmingham, ymode = log] 
			\addplot table[x index = 0, y index = 1] {make_birmingham.dat};
			\addplot table[x index = 0, y index = 2] {make_birmingham.dat};
			\addplot table[x index = 0, y index = 3] {make_birmingham.dat};
			\addplot table[x index = 0, y index = 4] {make_birmingham.dat};
			\addplot[mark = triangle, green] table[x index = 0, y index = 5] {make_birmingham.dat};
			\nextgroupplot[title=ChicagoRegional, ymode = log]
			\addplot table[x index = 0, y index = 1] {make_chicago.dat};
			\addplot table[x index = 0, y index = 2] {make_chicago.dat};
			\addplot table[x index = 0, y index = 3] {make_chicago.dat};
			\addplot table[x index = 0, y index = 4] {make_chicago.dat};
			\addplot[mark = triangle, green] table[x index = 0, y index = 5] {make_chicago.dat};
			\nextgroupplot[title=Hawaii, ylabel= Magnitude of the relative trace variation,	every axis y label/.append style={at=(ticklabel cs:1.1)}, xlabel = Budget of edges to be added ($k$), every axis x label/.append style={at=(ticklabel cs:1.75)}, ymode = log]
			\addplot table[x index = 0, y index = 1] {make_hawaii.dat};
			\addplot table[x index = 0, y index = 2] {make_hawaii.dat};
			\addplot table[x index = 0, y index = 3] {make_hawaii.dat};
			\addplot table[x index = 0, y index = 4] {make_hawaii.dat};
			\addplot[mark = triangle, green] table[x index = 0, y index = 5] {make_hawaii.dat};
			\nextgroupplot[title=RhodeIsland, ymode = log]
			\addplot table[x index = 0, y index = 1] {make_rhodeisland.dat};
			\addplot table[x index = 0, y index = 2] {make_rhodeisland.dat};
			\addplot table[x index = 0, y index = 3] {make_rhodeisland.dat};
			\addplot table[x index = 0, y index = 4] {make_rhodeisland.dat};
			\addplot[mark = triangle, green] table[x index = 0, y index = 5] {make_rhodeisland.dat};
			
			\nextgroupplot[title=Rome, ymode = log]
			\addplot table[x index = 0, y index = 1] {make_rome.dat};
			\addplot table[x index = 0, y index = 2] {make_rome.dat};
			\addplot table[x index = 0, y index = 3] {make_rome.dat};
			\addplot table[x index = 0, y index = 4] {make_rome.dat};
			\addplot[mark = triangle, green] table[x index = 0, y index = 5] {make_rome.dat};
		\end{groupplot}
	\end{tikzpicture}
	\caption{Magnitude of the relative trace variation for addition as the budget increases.}\label{fig:make_score}
\end{figure}

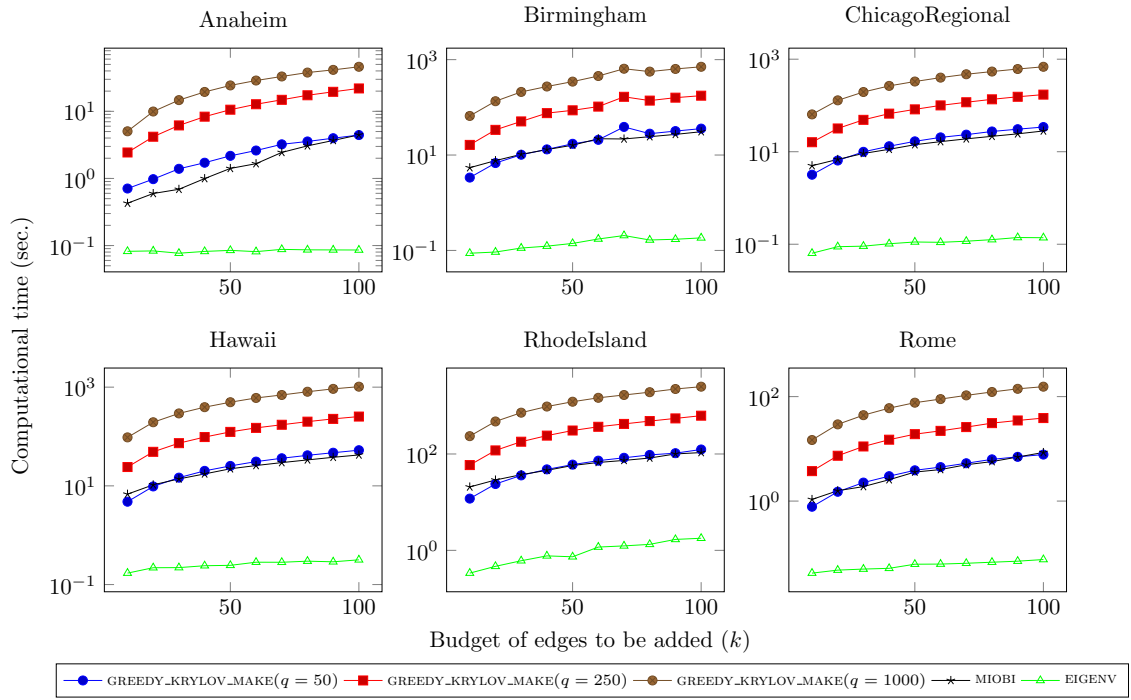
\begin{figure}
	\centering
	\begin{tikzpicture}[scale = 0.85]
	\begin{groupplot}[group style={group size= 3 by 2, vertical sep=1.5cm}, legend style={					
			font=\scriptsize,
			at={(3.69,-1.75)},
			legend columns = -1,
			cells={anchor=west}
		}, width=.4\textwidth]
		\nextgroupplot[title=Anaheim, scaled y ticks=base 10:1, ymode=log]
		\addplot table[x index = 0, y index = 1] {make_anaheim_time.dat};\addlegendentry{\textsc{greedy\_krylov\_make}($q=50$)};
		\addplot table[x index = 0, y index = 2] {make_anaheim_time.dat};\addlegendentry{\textsc{greedy\_krylov\_make}($q=250$)};
		\addplot table[x index = 0, y index = 3] {make_anaheim_time.dat};\addlegendentry{\textsc{greedy\_krylov\_make}($q=1000$)};
		\addplot table[x index = 0, y index = 4] {make_anaheim_time.dat};\addlegendentry{\textsc{miobi}};
		\addplot[mark = triangle, green] table[x index = 0, y index = 5] {make_anaheim_time.dat};\addlegendentry{\textsc{eigenv}};
		\nextgroupplot[title=Birmingham, ymode=log]
		\addplot table[x index = 0, y index = 1] {make_birmingham_time.dat};
		\addplot table[x index = 0, y index = 2] {make_birmingham_time.dat};
		\addplot table[x index = 0, y index = 3] {make_birmingham_time.dat};
		\addplot table[x index = 0, y index = 4] {make_birmingham_time.dat};
		\addplot[mark = triangle, green] table[x index = 0, y index = 5] {make_birmingham_time.dat};
		\nextgroupplot[title=ChicagoRegional, ymode=log]
		\addplot table[x index = 0, y index = 1] {make_chicago_time.dat};
		\addplot table[x index = 0, y index = 2] {make_chicago_time.dat};
		\addplot table[x index = 0, y index = 3] {make_chicago_time.dat};
		\addplot table[x index = 0, y index = 4] {make_chicago_time.dat};
		\addplot[mark = triangle, green] table[x index = 0, y index = 5] {make_chicago_time.dat};
		\nextgroupplot[title=Hawaii, ylabel= Computational time \upd{(sec.)},	every axis y label/.append style={at=(ticklabel cs:1.1)}, xlabel = Budget of edges to be added ($k$), every axis x label/.append style={at=(ticklabel cs:1.75)}, ymode=log]
		\addplot table[x index = 0, y index = 1] {make_hawaii_time.dat};
		\addplot table[x index = 0, y index = 2] {make_hawaii_time.dat};
		\addplot table[x index = 0, y index = 3] {make_hawaii_time.dat};
		\addplot table[x index = 0, y index = 4] {make_hawaii_time.dat};
		\addplot[mark = triangle, green] table[x index = 0, y index = 5] {make_hawaii_time.dat};
		\nextgroupplot[title=RhodeIsland, ymode=log]
		\addplot table[x index = 0, y index = 1] {make_rhodeisland_time.dat};
		\addplot table[x index = 0, y index = 2] {make_rhodeisland_time.dat};
		\addplot table[x index = 0, y index = 3] {make_rhodeisland_time.dat};
		\addplot table[x index = 0, y index = 4] {make_rhodeisland_time.dat};
		\addplot[mark = triangle, green] table[x index = 0, y index = 5] {break_rhodeisland_time.dat};
		
		\nextgroupplot[title=Rome, ymode=log]
		\addplot table[x index = 0, y index = 1] {make_rome_time.dat};
		\addplot table[x index = 0, y index = 2] {make_rome_time.dat};
		\addplot table[x index = 0, y index = 3] {make_rome_time.dat};
		\addplot table[x index = 0, y index = 4] {make_rome_time.dat};
		\addplot[mark = triangle, green] table[x index = 0, y index = 5] {make_rome_time.dat};
	\end{groupplot}
\end{tikzpicture}
	\caption{Computational times \upd{(seconds)} of the methods for addition as the budget increases.}\label{fig:make_time}
\end{figure}

Then, we investigate the impact of varying the budget size $k$ in the range $10,20,\dots,100$ on the trace variation and the computational time for the road networks considered in section~\ref{sec:exp-budget}. Also in this case, we consider three different sizes for the search space of \textsc{greedy\_krylov\_make}, corresponding to the choices of the parameter $q$ in the set of values $50,250,\min\{1000, |E|\}$. Figure~\ref{fig:make_score} reports the magnitude of the relative trace variation and highlights a crucial difference with respect to the downgrading problem: For any size of the search space, \textsc{greedy\_krylov\_make} outperforms significantly its competitors on all case studies. We also note that, in contrast to the downgrading case, \textsc{eigenv} has either comparable or better performances than \textsc{miobi} on all case studies.  Moreover, the computational times reported in Figure~\ref{fig:make_time} demonstrate that by choosing the smallest size of the search space ($q=50$), the cost of \textsc{greedy\_krylov\_make} becomes comparable to the one of \textsc{miobi}. This is also due to the fact that, for the addition problem, the search space of \textsc{miobi} might be significantly larger than in the downgrading case. Therefore, we conclude that \textsc{greedy\_krylov\_make} should be the method of choice for problem \eqref{eq:P-addition-binary}, unless a very strict limitation on the time consumption has to be applied.

\begin{table}[t]
    \centering
    \textbf{Addition}
    \resizebox{\textwidth}{!}{\begin{tabular}{lc@{\hskip .5em}c@{\hskip .5em}c@{\hskip .5em}c c@{\hskip 2em} lc@{\hskip .5em}c@{\hskip .5em}c@{\hskip .5em}c}
    \toprule
    & GKM&MIOBI&EIGENV&$\cap$  &&&GKM&MIOBI&EIGENV&$\cap$\\
    \cmidrule{1-5}  \cmidrule{7-11}
    Anaheim&\textbf{42.4}&12.8&15.9&31&&Cardiff&\textbf{24.0}&20.6&20.6&37\\
    Austin&\textbf{3.49}&2.34&2.34&37&&CollegeMsg&\textbf{5.14}&5.04&5.04&43\\
    Barcelona&\textbf{29.5}&12.3&12.3&30&&Edinburgh&\textbf{54.2}&19.9&19.9&24\\
    Birmingham&\textbf{1.36}&0.353&0.372&19&&as735&1.12&\textbf{2.31}&\textbf{2.31}&10\\
    ChicagoRegional&\textbf{1.64}&0.558&0.613&21&&AstroPh&\textbf{1.28}&1.27&1.27&39\\
    DC&\textbf{2.13}&0.197&0.495&19&&CondMat&4.66&\textbf{4.94}&\textbf{4.94}&30\\
    Hawaii&\textbf{1.13}&0.0745&0.271&16&&HepTh&\textbf{3.16}&2.78&2.78&39\\
    Philadelphia&\textbf{1.43}&0.162&0.293&14&&London&\textbf{70.9}&26.9&26.9&36\\
    RhodeIsland&\textbf{0.469}&0.150&0.150&27&&netscience&\textbf{52.1}&30.5&30.5&25\\
    Rome&\textbf{6.33}&2.34&2.17&22&&Epinions1&1.04&\textbf{1.13}&\textbf{1.13}&28\\
    Sydney&\textbf{0.794}&0.274&0.274&38&&yeast&\textbf{23.7}&20.4&20.4&32\\
    \bottomrule
    \end{tabular}}
    \caption{Magnitude of the relative trace variation obtained with the three methods \textsc{greedy\_krylov\_make (gkm)}, \textsc{miobi}, \textsc{eigenv} considered for the addition of edges to unweighted graphs on road networks (left) and general networks (right), with a budget of $k=50$ edges. The column denoted with $\cap$ shows the number of edges that have been commonly chosen  by all the methods.}
    \label{tab:make_road_misc}
\end{table}

\section{Numerical experiments with weighted graphs: tuning, rewiring, addition}\label{sec:exp-weight}


Finally, we present results on a set of weighted networks in order to test the performance of the proposed method for the edge tuning problem \eqref{eq:P-tuning}, as well as weighted edge-addition and edge rewiring, where we simultaneously tune the weight of existing edges and add new ones. 

While in certain applications the set $F$ of edges (or missing edges) that we are allowed to modify  is given a-priori, in our setup we will assume only the cardinality of the set  $F$ is fixed, i.e.\ we are free to select a  set  of $n_F\geq 1$ modifiable edges (or edges to be added) and we need to form $F$ by choosing which ones are those that are best suited to maximize the natural connectivity. This is a more challenging scenario and, clearly, the case in which the set $F$ is specified by external constraints is retrieved as a special case.  

In order to find an optimal set $F$, we propose to measure how sensitive the \upd{$f$-connectivity, with $f\in\{\exp,\ \sinh\}$,} with respect to changes in the weight of a certain edge $(i,j)$ is. To this end, one should look at the magnitude of the corresponding gradient entry \upd{$2f'(A+X)_{ij}$} and  select the edges corresponding to the largest gradient.  However, inspecting these quantities for all edges (or missing edges) can be too expensive for large networks. Thus, in our experiments, we proceed as follows:
    first, we select a set of $n_P$ candidate edges, with $n_P > n_F$, chosen as the most important, with respect to a suitable edge-ordering, among existing and/or non-existing edges;
    then, we identify $F$ on the basis of the evaluations of the gradient over the $n_P$ candidate edges. 


\begin{table}[t!]
	\caption{Number of vertices and edges of the weighted graphs of power grids used for the numerical tests. In bracket the percentage of voltages that has been added in our preprocessing stage.}\label{tab:info-power-grids}
	\vspace{.2cm}
	
	\centering
	\begin{tabular}{lcc}
		\toprule	
		\textbf{Dataset}& $|V|$&$|E|$\\
		\midrule
		Austria&$149$&$169$ $(0\%)$\\
		Denmark&$96$&$105$ ($0\%$)\\
		England&$504$&$603$ ($0.7\%$)\\
		Germany&$1903$& $2371$ ($0.4\%$)\\
		Italy&$858$&$1092$ ($2.8\%$)\\
		India&$3228$&$4323$ ($0.4\%$)\\
		Mexico&$552$&$743$ ($2.8\%$)\\
		Poland&$299$&$390$ ($0.3\%$)\\
		Portugal&$185$&$247$ $(1.5\%)$\\
		Sweden&$268$&$336$ ($2.4\%$)\\
		\bottomrule
	\end{tabular}
\end{table}

In our experiments, we test \textsc{krylov\_lbfgs} \upd{and \textsc{krylov\_hessian}} on a set of electric power grid networks from different countries, 
as listed in Table~\ref{tab:info-power-grids}. All the considered network datasets were collected from an Open Street Map project by the Complex Network Group at Telecom Sud-Paris  \cite{power_grid_data}. Each node represents a power station and edges represent wired connections, weighted by their voltage capacity. A small number (in most cases less than 1\%) of edge voltage capacity data was missing in the original datasets. For those edges we artificially set the voltage capacity as the average of the neighbors. \upd{In all the tests of this section, we consider the total weight budget $k=10$.}

Concerning the selection of the edges in $F$, we propose three different approaches that deal with different scenarios, as listed below.

\begin{paragraph}{Tuning.}
This approach applies to the case where we are only allowed to modify edges with an initial non-zero weight. We select the candidate edges as the first $n_P=100$ existing edges with respect to $\le_1$; then, we set $F$ as the $n_F=30$ edges, among the candidates, with the largest value of the gradient. 
\end{paragraph}  
\begin{paragraph}{Rewiring.}
This approach applies to the case where we are allowed  to both modify existing edges and add new ones. We select two sets $C_1$ and $C_2$ of $50$ candidate edges each, as the first $50$  existing edges with respect to $\le_2$ and the first $50$ non-existing edges with respect to $\le_2$, respectively. The resulting set of $n_P=100$ candidate pairs is then used to form $F$ by choosing  $n_F=30$ elements from the union of the $15$ edges in $C_1$ and non-existing edges in $C_2$  with the largest value of the gradient. 
\end{paragraph}
\begin{paragraph}{Addition.}
This approach applies to the case where we are only allowed to add new edges. We select the candidate edges as the first $n_P=100$ non-existing edges with respect to $\le_2$; then, we set $F$ as the $n_F=30$ edges, among the candidates, with the largest value of the gradient. 
\end{paragraph}
\vspace{.2cm}

\begin{figure}[t]
    \centering
    \includegraphics[width=\textwidth]{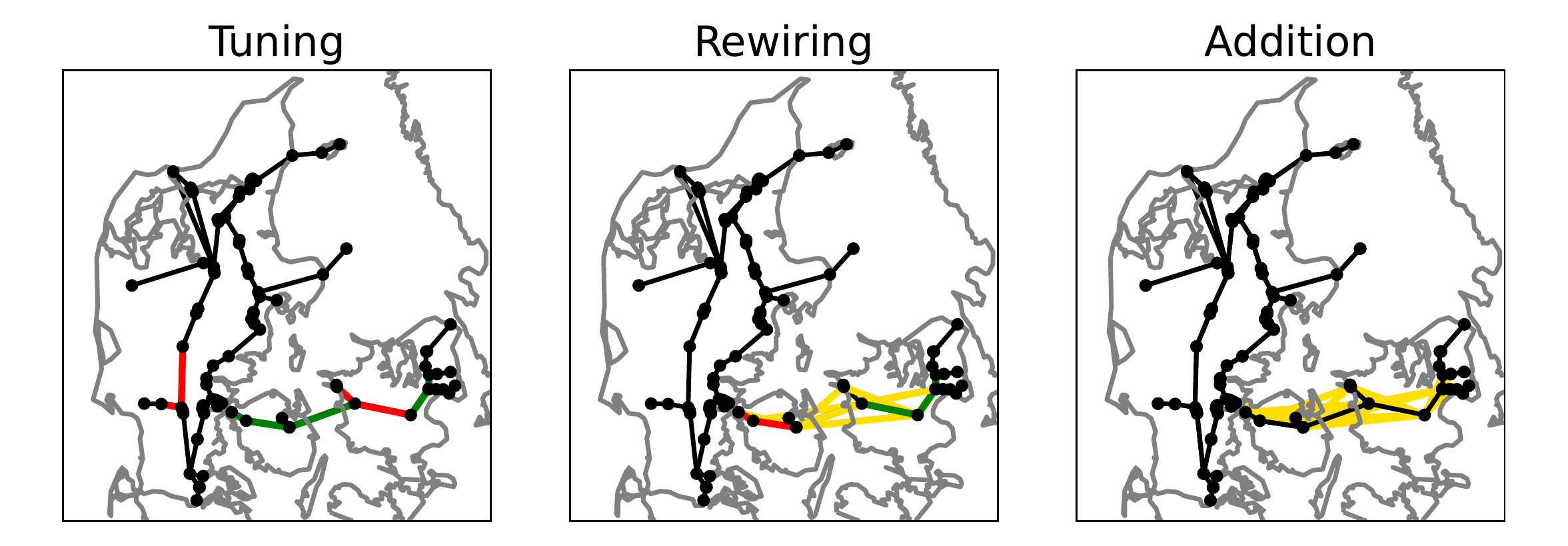}
    \caption{Results of \textsc{krylov\_lbfgs} on the electric power grid network of Denmark. Red lines correspond to edges whose weight was decreased; Green lines correspond to edges whose weight was increased; Yellow lines denote edges that have been added from scratch.}
    \label{fig:Denmark}
\end{figure}

\upd{Tables \ref{tab:tune_exp} shows the relative trace variation and the execution time (in seconds), \upd{obtained with \textsc{krylov\_lbfgs} and \textsc{krylov\_hessian}}, for all the datasets and the three problem cases above. The values of $\Delta T$ obtained with the two methods are very close, indeed their difference is more than the $10$\% of the highest value only in two cases: Austria (Rewiring) and Portugal (Addition).} \Cref{fig:Denmark} shows the geographical location of the modified and added edges on the power network of Denmark, obtained with \textsc{krylov\_lbfgs}, where red edges denote edges whose weight has been diminished by the algorithm, green edges are edges whose weight was increased, and yellow lines denote edges that were added. As expected, rewiring is always the most effective procedure, resulting in the largest increase in natural connectivity, as it combines edge tuning and edge addition in a simultaneous optimization mechanism. In particular, we see from \Cref{fig:Denmark} that the set of edges modified and added by Rewiring is a subset of those that are modified and added by the other two approaches. 
	
	\upd{Empirically, we observe that \textsc{krylov\_hessian} always converge in less iterations; however, the latter are more expensive and there is no clear winner between the two methods, in terms of speed; \textsc{krylov\_lbfgs} is faster on $16$ examples while \textsc{krylov\_hessian} on $14$.}

\upd{The numerical test is repeated with $f=\sinh$ and the corresponding results are reported in Table \ref{tab:tune_sinh}. On all case studies \textsc{krylov\_lbfgs} and \textsc{krylov\_hessian} yields almost equal variations of the $f$-connectivity. The most significant differences are observed on England (Addition), Germany (Rewiring), Mexico (Addition), and Poland (Tuning). Similar comments to the exponential case, apply to the reported computational times.} 

\begin{table}[t]
	\centering
	\textbf{Exponential}
	\resizebox{\textwidth}{!}{\begin{tabular}{l c@{\hskip .4em}c@{\hskip .4em}c   c@{\hskip .1em}   c@{\hskip  .4em}c@{\hskip  .4em}c   c@{\hskip  .5em} c@{\hskip  .4em}c@{\hskip  .4em}c     c@{\hskip  .1em} c@{\hskip  .4em}c@{\hskip  .4em}c    c@{\hskip  .5em} c@{\hskip  .4em}c@{\hskip  .4em}c c@{\hskip  .1em} c@{\hskip  .4em}c@{\hskip  .4em}c}
    \toprule
    & \multicolumn{7}{c}{\textsc{Tuning}} && \multicolumn{7}{c}{\textsc{Rewiring}} && \multicolumn{7}{c}{\textsc{Addition}} \\
    \cmidrule{2-8} \cmidrule{10-16} \cmidrule{18-24}
    & \multicolumn{3}{c}{\textsc{lbfgs}} && \multicolumn{3}{c}{\textsc{Hessian}} && \multicolumn{3}{c}{\textsc{lbfgs}} && \multicolumn{3}{c}{\textsc{Hessian}} && \multicolumn{3}{c}{\textsc{lbfgs}} && \multicolumn{3}{c}{\textsc{Hessian}} \\
    \cmidrule{2-4} \cmidrule{6-8} \cmidrule{10-12} \cmidrule{14-16} \cmidrule{18-20} \cmidrule{22-24}
    &$\Delta T$ & t & it &&$\Delta T$ & t & it &&$\Delta T$ & t & it &&$\Delta T$ & t & it && $\Delta T$ & t & it && $\Delta T$ & t & it  \\
    \midrule
    Austria&0.76&0.2&29&&\textbf{0.76}&3.2&15&&2.89&0.4&34&&\textbf{3.2}&1.3&12&&\textbf{0.55}&0.2&20&&\textbf{0.55}&1.5&11\\
    Denmark&\textbf{1.22}&0.2&47&&1.2&1.5&14&&\textbf{6.42}&0.1&34&&\textbf{6.42}&1.4&12&&\textbf{0.71}&0.2&28&&\textbf{0.71}&1.2&13\\
    England&\textbf{0.32}&3.2&34&&\textbf{0.32}&3.9&14&&0.8&9.9&54&&\textbf{0.8}&4.5&16&&0.2&1.3&19&&\textbf{0.22}&2.1&13\\
    Germany&0.19&6.5&42&&\textbf{0.19}&4.9&17&&\textbf{0.53}&2.4&45&&\textbf{0.53}&2.1&12&&\textbf{0.08}&1.2&22&&\textbf{0.08}&2.0&12\\
    India&\textbf{0.12}&9.0&38&&\textbf{0.12}&6.5&15&&\textbf{0.34}&7.5&53&&0.33&4.4&16&&\textbf{0.06}&2.3&17&&\textbf{0.06}&2.6&10\\
    Italy&\textbf{0.41}&7.3&37&&\textbf{0.41}&4.1&12&&\textbf{1.62}&4.6&46&&\textbf{1.62}&2.7&12&&\textbf{0.16}&2.3&27&&\textbf{0.16}&2.0&12\\
    Mexico&\textbf{0.58}&7.7&40&&\textbf{0.58}&4.2&12&&1.66&13.6&47&&\textbf{1.66}&4.5&12&&0.2&1.4&18&&\textbf{0.21}&2.0&11\\
    Poland&\textbf{0.62}&2.6&37&&\textbf{0.62}&2.1&13&&\textbf{1.76}&1.8&34&&\textbf{1.76}&1.7&10&&\textbf{0.31}&1.8&31&&\textbf{0.31}&2.1&14\\
    Portugal&\textbf{1.02}&0.5&42&&\textbf{1.02}&1.6&14&&\textbf{3.46}&0.5&27&&\textbf{3.46}&1.4&12&&\textbf{0.6}&0.6&35&&0.47&1.8&19\\
    Sweden&0.6&1.0&33&&\textbf{0.62}&1.6&13&&2.49&0.8&26&&\textbf{2.49}&1.4&12&&0.33&1.2&33&&\textbf{0.33}&1.5&13\\
    \bottomrule
    \end{tabular}}
	\caption{Magnitude of the relative trace variation ($\Delta T$), execution time in seconds (t), and number of iterations (it),  for the tuning, rewiring, and addition optimization problems solved with \textsc{krylov\_lbfgs} and \textsc{krylov\_hessian} approaches, for weighted graphs associated with power grid networks and for $f=\exp$.}
	\label{tab:tune_exp}
\end{table}

\begin{table}[t]
    \centering
    \textbf{Hyperbolic sine}
    \resizebox{\textwidth}{!}{\begin{tabular}{l c@{\hskip .4em}c@{\hskip .4em}c   c@{\hskip .1em}   c@{\hskip  .4em}c@{\hskip  .4em}c   c@{\hskip  .5em} c@{\hskip  .4em}c@{\hskip  .4em}c     c@{\hskip  .1em} c@{\hskip  .4em}c@{\hskip  .4em}c    c@{\hskip  .5em} c@{\hskip  .4em}c@{\hskip  .4em}c c@{\hskip  .1em} c@{\hskip  .4em}c@{\hskip  .4em}c}
    \toprule
    & \multicolumn{7}{c}{\textsc{Tuning}} && \multicolumn{7}{c}{\textsc{Rewiring}} && \multicolumn{7}{c}{\textsc{Addition}} \\
    \cmidrule{2-8} \cmidrule{10-16} \cmidrule{18-24}
    & \multicolumn{3}{c}{\textsc{lbfgs}} && \multicolumn{3}{c}{\textsc{Hessian}} && \multicolumn{3}{c}{\textsc{lbfgs}} && \multicolumn{3}{c}{\textsc{Hessian}} && \multicolumn{3}{c}{\textsc{lbfgs}} && \multicolumn{3}{c}{\textsc{Hessian}} \\
    \cmidrule{2-4} \cmidrule{6-8} \cmidrule{10-12} \cmidrule{14-16} \cmidrule{18-20} \cmidrule{22-24}
    &$\Delta T$ & t & it &&$\Delta T$ & t & it &&$\Delta T$ & t & it &&$\Delta T$ & t & it && $\Delta T$ & t & it && $\Delta T$ & t & it  \\
    \midrule
    Austria&\textbf{13.98}&0.9&41&&13.9&1.4&9&&\textbf{67.14}&0.4&26&&67.14&1.4&10&&19.85&0.3&15&&\textbf{19.85}&1.0&9\\
    Denmark&\textbf{15.51}&0.3&38&&15.5&1.1&11&&123.74&0.2&23&&\textbf{123.75}&1.3&11&&16.72&0.1&19&&\textbf{16.74}&1.0&8\\
    England&\textbf{3.68}&10.0&40&&3.68&6.0&13&&\textbf{9.7}&4.2&35&&9.7&2.9&10&&3.29&2.0&20&&\textbf{3.67}&2.2&10\\
    Germany&2.72&3.0&12&&\textbf{2.74}&2.1&7&&9.26&2.6&39&&\textbf{9.87}&2.1&11&&1.19&1.9&24&&\textbf{1.32}&1.8&10\\
    India&\textbf{6.86}&6.8&24&&6.86&4.9&12&&\textbf{32.21}&5.0&29&&32.17&3.2&9&&5.26&3.5&24&&\textbf{5.26}&2.7&9\\
    Italy&5.48&11.1&47&&\textbf{5.54}&4.1&12&&19.46&3.2&20&&\textbf{19.46}&3.2&10&&2.8&2.1&30&&\textbf{2.8}&2.4&10\\
    Mexico&\textbf{5.02}&4.9&25&&5.01&3.0&9&&\textbf{12.75}&7.4&33&&12.72&3.1&9&&2.05&3.6&35&&\textbf{2.6}&3.2&15\\
    Poland&6.23&1.4&17&&\textbf{6.44}&1.7&9&&14.45&2.7&28&&\textbf{14.46}&2.2&10&&\textbf{4.05}&3.2&27&&4.0&1.5&7\\
    Portugal&9.67&1.1&27&&\textbf{9.68}&2.6&10&&40.63&0.8&24&&\textbf{41.3}&1.3&10&&\textbf{6.41}&0.9&29&&6.4&1.5&12\\
    Sweden&\textbf{6.9}&1.7&27&&6.9&7.2&34&&\textbf{43.36}&1.8&24&&43.09&1.6&8&&\textbf{5.27}&1.7&25&&5.27&1.5&8\\
    \bottomrule
    \end{tabular}}
    \caption{Magnitude of the relative trace variation ($\Delta T$), execution time in seconds (t), and number of iterations (it), for the tuning, rewiring, and addition optimization problems solved with \textsc{krylov\_lbfgs} and \textsc{krylov\_hessian} approaches, for weighted graphs associated with power grid networks and for $f=\sinh$.}
    \label{tab:tune_sinh}
\end{table}

\section{Conclusions}
We have proposed two strategies, based on Krylov subspace approximations, for optimizing the natural connectivity of a graph. The first one is a greedy heuristic method that is well suited to contexts where we either add or remove unweighted edges on a large-scale graph. \upd{Despite been computationally more expensive than state-of-the-art alternatives, in the context of the addition problem our approach significantly outperforms the increase of the natural connectivity}. The second proposed strategy combines Krylov subspace approximation and \upd{an interior point scheme using either the Hessian or its L-BFGS approximation,} to address continuous optimization problems that include edge tuning and rewiring. To the best of our knowledge, this is the first attempt to tackle the optimization of the natural connectivity with \upd{first and second order methods}, and the reported experiments demonstrate the feasibility of the approach at least for graphs up to medium size. 

Finally, we highlight that the proposed computational strategies are quite flexible as they can be adapted with minor changes to the optimization of other matrix function based measures on graphs and it is conceptually easy to incorporate further constraints on the set of modifiable edges.

\section*{Acknowledgments}
We would like to thank the department of Math and Stats of Uni Strathclyde for hosting us and the  European Union's Horizon 2020 research and innovation programme who has provided support for the researchers to engage in collaborating activities via the   Marie Sk\l odowska-Curie individual fellowship ``MAGNET'' No 744014.

\bibliographystyle{abbrv} 
\bibliography{references}  

\end{document}